\documentclass[a4paper]{article}
\usepackage[english]{babel}
\usepackage{graphicx}

\newcommand \ve  \varepsilon

\def\M{{\mathcal M}}
\def\S{{\mathcal S}}
\def\iff{\textrm{if}\hskip 2mm }
\def\N{N^*}
\def\NN{{\mathcal N}}
\def\Z{{\mathcal Z}}
\def\ZZ{{\mathcal Z}^2}
\def\I{{\mathcal J}}
\def\E{{\mathcal E}}
\def\A{{\mathcal A}}
\def\L{{\mathcal L}}

\def\bx{{\bf x}}
\def\bt{{\bf \theta }}

\begin{document}

\begin{center}
\Large

{\bf Explicit Construction of }

{\bf the Brownian Self-Transport Operator }

\vskip 3mm
D.S.Grebenkov\footnote{ E-mail: Denis.Grebenkov@polytechnique.fr }
\normalsize
\vskip 5mm
Laboratoire de Physique de la Mati\`ere Condens\'ee, Ecole Polytechnique

91128 Palaiseau Cedex \hskip 1mm France
\vskip 3mm

Department of Statistical Physics, Saint Petersburg State University

ul. Ulyanovskaya 1, Petrodvorets, 198904, Saint Petersburg, Russia

\end{center}

\vskip 10mm
\begin{abstract}

Applying the technique of characteristic functions developped for
one-dimensional regular surfaces (curves) with compact support, we 
obtain the distribution of hitting probabilities for a wide class of 
finite membranes on square lattice. Then we generalize it to 
multi-dimensional finite membranes on hypercubic lattice. 
Basing on these distributions, we explicitly construct 
the Brownian self-transport operator which governs the Laplacian transfer. 

In order to verify the accuracy of the distribution
of hitting probabilities, numerical analysis is carried out 
for some particular membranes.

\end{abstract}

\vskip 10mm
\section*{ Introduction }

We are engaged in the physical problems that can be mathematically 
modelized by {\it the Laplacian transfer problem} \cite{Sapoval94}. 
These are : the diffusion through semi-permeable membranes and the electrode 
problem \cite{Sapoval99}, the heterogeneous catalysis \cite{Sapoval01}, 
the NMR in porous environment \cite{Sapoval96}, etc. 

\begin{figure}             
\begin{center}
\includegraphics[width=12cm,height=45mm]{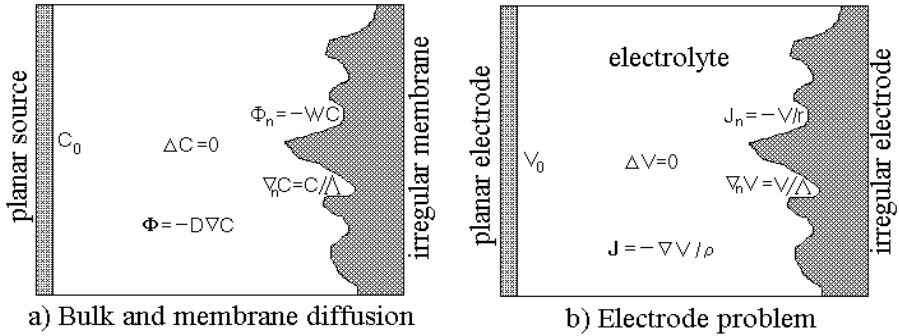}   
\end{center}
\caption{ The problem of the Laplacian transfer across a resistive
irregular interface (a membrane or an electrode). }                             \label{fig:problem}
\end{figure}

Consider the diffusion of particles from the source (usually supposed planar) 
to the semi-permeable membrane with rather complex geometry, see Fig.~\ref{fig:problem}a. 
In the steady-state regime, the concentration of particles $C$ obeys the 
Laplacian equation,
\begin{equation}                                                                 \label{eq:Laplace}
\Delta C=0 .
\end{equation}
The flux of particles in the bulk obeys Fick's law, ${\bf \Phi}=-D{\bf \nabla }C$,
where $D$ is the diffusion coefficient. The flux across the surface is
given by $\Phi _n=-WC$, where $W$ is the permeability of the membrane 
(the probability per unit time, surface, and concentration for a particle
to cross the membrane). Equating these two fluxes, we obtain 
the mixed boundary condition,
\begin{equation}                                                                 \label{eq:boundary}
\frac{\partial C}{\partial n}=\frac{1}{\Lambda }C ,
\end{equation}
often called also Fourier or Robin boundary condition. The physical parameter
$\Lambda =D/W$ plays an important role in this treatment.

In the electrode problem (see Fig.~\ref{fig:problem}b), there are a planar and an irregular 
electrodes (with an interface resistivity $r$), placed in an electrolyte 
with a resistivity $\rho$. The electric potential $V$ obeys the Laplace 
equation $\Delta V=0$ in the bulk of the electrolyte. The boundary condition 
(\ref{eq:boundary}) is obtained by equating the current from the electrolyte, 
$-\nabla V/\rho $, to the current $-V/r$ crossing the electrode surface, 
with $\Lambda =r/\rho $.

For the heterogeneous catalysis, a catalyst of a complex geometry
(with the reactivity $K$) is placed in a solution. When a molecule $A$
diffusing in the solution, reaches the catalytic surface, it transforms 
into $A^*$ with reaction rate $K$. As above, the concentration of molecules $A$
obeys the Laplace equation (\ref{eq:Laplace}). The boundary condition 
(\ref{eq:boundary}) with $\Lambda =D/K$ is obtained due to the mass conservation.

The same arguments allow to describe the NMR in the porous environment
by the Laplace equation under the mixed boundary condition. 
So, all these different phenomena are described by the same mathematical 
formalism (see \cite{Filoche99}, \cite{Grebenkov00} for details). 

\vskip 2mm

The Laplacian transfer problem is significantly more complex 
than the corresponding Dirichlet or Neumann problem.
In order to find the solution, Filoche and Sapoval \cite{Filoche99} 
proposed the following program :

  -- to choose the appropriate discretization;

  -- to solve the discrete problem;
  
  -- to take the continuous limit.

For these purposes, they introduced {\it the Brownian self-transport operator} $Q$.
This operator controls the properties of the Laplacian 
transfer and depends only on the geometry of membrane's boundary.
It can be easily defined for a given discrete membrane: 
$Q_{k,n}$ is the probability of first contact with the membrane on the $n$-th 
site if started from the $k$-th site (without touching other sites!). 
This is exactly the problem of the hitting probabilities on the lattice. 
So, if one can calculate $Q$ for any ``reasonable'' membrane, the Laplacian 
transfer problem would be solved (except the difficulties related with 
the continuous limit). In this paper, we construct the Brownian self-transport
operator for the membranes with a rather general geometry.

\begin{figure}             
\begin{center}
\includegraphics[width=11cm,height=3cm]{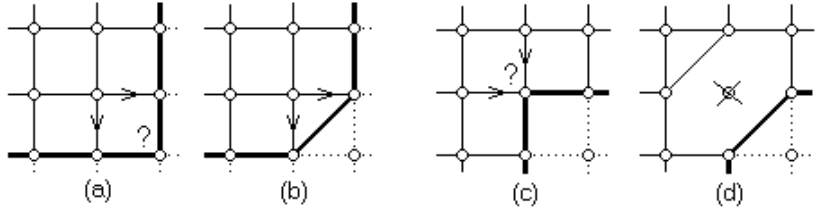}
\end{center}
\caption{ Corners: (a) internal corner;
(b) corner point is carried from boundary to the bulk of membrane;
(c) external corner; (d) corner point is removed, two links
from boundary to corner point are also removed, two other links
are connected to each other.}                                                    \label{fig:corner}
\end{figure}

Now we should specify what the term ``reasonable membrane'' means.
For the moment, we confine ourselves to the two-dimensional case 
(square lattice) in order to clarify the representation of our 
treatment. When two-dimensional results are obtained,
we shall explain how they can be generalized for the multi-dimensional
case. On a square lattice, the Laplacian equation (\ref{eq:Laplace}) 
becomes
\begin{equation}                                                        \label{eq:Laplace_discrete}
4C_{x,y}-C_{x+1,y}-C_{x-1,y}-C_{x,y+1}-C_{x,y-1}=0,
\end{equation}
i.e., diffusion can be modelized by {\em simple} random walks. 
We suppose that the membranes are quite regular. 
It means that one can choose such discretization that the
number of singular points is negligible by comparison with
the number of regular points (contrariwise to the case of 
fractal membranes). As we decided to work on the square lattice,
we should represent the membrane as a sequence of horizontal
and vertical segments. The last difficulty concerns with
corners. The problem is that a corner point is connected with 
{\em two} lattice points whereas a regular point is  
connected with only one lattice point. We would like to 
be indemnified from these ambiguities\footnote{ The other reason
to avoid the corner points is related with numerical simulations.
Indeed, to verify our analytical results we should solve the problem
of Laplacian transfer numerically. And here it is not clear
how to treat the corner points.}. To avoid the corner points,
we use the following rules\footnote{ These rules were described,
in particular, in \cite{Grebenkov00}. }:

-- for the internal corner (see Fig.~\ref{fig:corner}a,~\ref{fig:corner}b), 
we carry the corner point from the boundary to the bulk of membrane.

-- for the external corner (see Fig.~\ref{fig:corner}c,~\ref{fig:corner}d), 
we completely remove the corner point from the lattice. The links from the boundary
to these points are also removed. Two other links are closed
to each other.

This operation preserves the connectivity of the lattice and 
eliminates the singularities (corner points).

Section~\ref{sec:definitions} introduces definitions and restrictions 
which allow us to reach final results. Section~\ref{sec:technique} is devoted to 
the technique developped in our previous paper \cite{Grebenkov02} for 
{\em one-dimensional regular surfaces with compact support}.
Note that the actual treatment is partially based on the previous results. 
The distribution of hitting probabilities for a two-dimensional finite membrane 
is obtained in Section~\ref{sec:2D}. The generalization of these results for 
multi-dimensional case is described in Section~\ref{sec:multi}. In Section~\ref{sec:BSO}, 
we explicitly construct the Brownian self-transport operator. Some important 
generalizations are given in Section~\ref{sec:generalizations}. Section~\ref{sec:numerical} 
describes some numerical results in order to verify the method. 
In the last section we make conclusions.

\vskip 10mm
\section{  Definitions }                                                     \label{sec:definitions}

Consider a square lattice $\ZZ$ on a plane. A {\it path} is a sequence 
(finite or infinite) of points $\{ A_n\}\subset \ZZ$
such that $\forall n\hskip 2mm |A_n-A_{n-1}|=1$, where $|A-B|$ denotes the
distance between two points $A$ and $B$.

The sequence of points $\S =\{ (x_m,y_m)\in \ZZ \: :\: m\in \Z \}$ is
called {\it membrane's boundary} if it obeys the following conditions : 

1. {\em Bijection}: Index $m$ enumerates the points $(x_m,y_m)$, i.e.,
there is one-to-one correspondance (bijection) between the set of 
integer numbers $\Z$ and the set of boundary points $\S$.

2. {\em Boundary}: $\S $ separates the lattice $\ZZ$ in
two disjoint\footnote{ I.e., there is no path from one set to
another which does not pass through a point in $\S $.} sets -- 
$\E$ and $\M _0$ -- {\it external} and {\it internal} points. 

3. {\em Accessibility}: Any point of $\S$ is accessible
from all external and internal points, i.e., for any $(x,y)\notin \S $
and any $(x_m,y_m)\in \S $ there exists a path $\{ A_n\}$ such that
$A_1=(x,y)$ and $\{ A_n\}\cap \S =(x_m,y_m)$.
 
4. {\em Compactness}: There is only {\it finite} number $M$ of 
boundary points (in $\S $) which do not lie on the horizontal axis, 
i.e., the non-plane part of the membrane has a finite size. 

Set $\M =\M _0\cup \S$ is called {\it finite membrane}.

Once chosen, points $(x_m,y_m)$ are assumed fixed
in all following calculations. We always use 
notations $x_m$ (or $x_n$) and $y_m$ (or $y_n$) for 
abscissae and ordinates of boundary points. 

\vskip 1mm

In order to simplify expressions and to avoid possible ambiguities, 
we introduce the following conventions:

-- the bulk of membrane is placed in the lower half plane (except,
possibly, a finite number of points);

-- non-plane boundary points are enumerated by index $m=1..M$, i.e.,
$$\forall \: m\in [1..M]\hskip 4mm  y_m\ne 0, \hskip 8mm 
\forall \: m\notin [1..M]  \hskip 4mm y_m=0  . $$

\begin{figure}              
\begin{center}
\includegraphics[width=108mm,height=3cm]{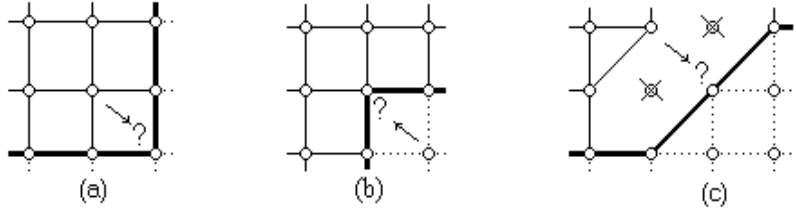}
\end{center}
\caption{ The accessibility condition forbids the corner points
(Fig.~\ref{fig:accessibility}a,~\ref{fig:accessibility}b) and 
``diagonal'' points (Fig.~\ref{fig:accessibility}c).}                       \label{fig:accessibility}
\end{figure}

Let us discuss the definition of the finite membrane.
The boundary and accessibility conditions provide that 
boundary points take their places sufficiently close
to each other, i.e., $\forall m\: \exists n\: |(x_m,y_m)-(x_n,y_n)|\leq \sqrt{2}$.
The accessibility condition prohibits the existence of corner points
(see Fig.~\ref{fig:accessibility}a,~\ref{fig:accessibility}b). Moreover, 
this condition forbids also the ``diagonal'' points 
(see Fig.~\ref{fig:accessibility}c). This means that membrane's boundary 
is composed of horizontal and vertical segments
just as required. So, we can conclude that these two conditions
provide all the necessary properties which were described
in the Introduction. 

The compactness condition is essential to achieve our goal.
This condition means that the membrane is composed by three
parts: two plane ``tails'' on the same height, and one 
intermediate (non-trivial) part that can be very complex but finite.
The importance of the compactness condition was throughly discussed 
in the previous paper, and we do not repeat it here. Note that in practice, 
this condition is not too restrictive because usually one considers the 
finite number of sites. 

According to the bijection condition, we can define the function
$\I (x,y)$ which gives the index of the boundary points $(x,y)$, i.e.,
$\I (x_m,y_m)=m$. Note that this function is defined only on the
boundary points: if $(x,y)\notin \S$, $\I (x,y)$ has no value.
We shall use the following convention: 

\begin{center}
if $(x,y)\notin \S$, the object containing $\I (x,y)$ is equal to $0$. 
\end{center}
For example, if we write $e^{i\I (x,y)\theta }$, we always mean 
$e^{i\I (x,y)\theta }\chi _{\S}(x,y)$, where 
$$\chi _S(x,y)=\cases{ 1, \: \iff (x,y)\in \S ,\cr
                       0, \: \iff (x,y)\notin \S. }$$
Note that it is just a useful convention to simplify the expressions.

We call all the points $\{ (x,y)\in \ZZ\: :\: y=n\}$ the {\it $n^{th}$ level}. 
We say that the membrane's boundary lies between $(-\N)^{th}$ and $N^{th}$ levels if
$$N=\max \{ y_m \}, \hskip 15mm \N=-\min \{ y_m \} .$$

The external point $A$ is called {\it near-boundary point},
if there exists $m$\ such that $|A-(x_m,y_m)|=1$, i.e., $A$ lies ``near'' the 
membrane's bourdary. The functions defined on these points, are 
called {\it near-boundary functions} (see below). We enumerate 
these points (and functions) by the same index $m$ as for
boundary points. Note that usually we suppose $m\in [1,M]$, i.e.,
we take into account near-boundary points which lie near the non-plane part 
of the membrane.

The external points $(x,0)\notin \M$ are called 
{\it ground points}. The functions defined on these points, 
are called {\it ground functions}. Let us enumerate ground
points $(\bar{x}_g,0)$ by index $g\in [1,G]$ using bar notation 
for their abscissae. Note that there exist membranes without 
ground points (see Fig.~\ref{fig:2D}a for example). 

We introduce the {\it outer normal} $(\delta x_m,\delta y_m)$ for 
each site of the membrane's boundary. Due to the accessibility condition, 
we have no corner points, therefore this vector is correctly defined. 
Note that $m$-th near-boundary point is just $(x_m+\delta x_m,y_m+\delta y_m)$, 
i.e., the outer normal is a vector directed 
from the boundary point to the corresponding near-boundary point. 

Let $P_{x,y}(n)$ be the hitting probability, i.e., the probability of the first 
contact with the membrane on point $(x_n,y_n)$ if started from an external point 
$(x,y)\in \E$ without touching other boundary points. Their characteristic functions 
$\phi _{x,y}(\theta )$ are
$$\phi _{x,y}(\theta )=\sum\limits _{m=-\infty }^{\infty }P_{x,y}(m)
e^{im\theta }  .$$
The inverse Fourier transform allows to obtain $P_{x,y}(m)$,
$$P_{x,y}(m)=\int\limits _{-\pi }^{\pi }\frac{d\theta }{2\pi }e^{-im\theta }
\phi _{x,y}(\theta ) .$$

\vskip 10mm
\section{ Technique developped in \cite{Grebenkov02} }                        \label{sec:technique}
\vskip 2mm

Here we briefly present the technique which was developped in
\cite{Grebenkov02} to obtain the distribution of hitting
probabilities for a regular surface with compact support. 
All the details can be found in this paper. 

The discrete Laplacian equation (\ref{eq:Laplace_discrete}) in terms 
of characteristic functions is simply
\begin{equation}                                                                 \label{eq:relation}
\phi _{x,y}(\theta )=\frac14\biggl[\phi _{x+1,y}(\theta )+\phi _{x-1,y}(\theta )+
\phi _{x,y-1}(\theta )+\phi _{x,y+1}(\theta )\biggr] .
\end{equation}
Using the convolution properties of hitting probabilities and their normalization, 
one obtains the distribution of hitting probabilities for a planar surface,
\begin{equation}                                                                       \label{eq:H}
P^{(planar)}_{x,y}(m)=\int\limits _{-\pi }^{\pi }\frac{d\theta }{2\pi }
e^{i(x-x_m)\theta }\varphi ^{|y|}(\theta )=H^y_{x-x_m}, 
\end{equation}
where function $\varphi (\theta )$ is
\begin{equation}                                                                  \label{eq:varphi}
\varphi (\theta )=2-\cos \theta -\sqrt{(2-\cos \theta )^2-1}.
\end{equation} 
We wrote $\varphi ^{|y|}(\theta )$ given that for the lower half plane ($y<0$) 
one has exactly the same result due to reflection symmetry.

\vskip 1mm

For the non-planar membrane, relations (\ref{eq:relation}) remain valid for 
{\it external points}, but there are also {\it membrane's points} where it is not true. 
Nevertheless, we preserve this form by adding the correction term which
equals to $0$ on external points.

In order to simplify manipulation with characteristic functions, 
we introduce vectors $\Phi ^{(y)}(\theta )$ containing all $\phi _{x,y}(\theta )$ 
on the $y$-th level with $|x|\leq L$ (at the end of calculations $L\to \infty $). 
We can write (\ref{eq:relation}) as
$$A\Phi ^{(y)}(\theta )=\Phi ^{(y-1)}(\theta )+\Phi ^{(y+1)}(\theta )+
\Delta \Phi ^{(y)}(\theta ),$$
where $A$ is tridiagonal matrix with elements: $A_{i,i}=4$, $A_{i,i+1}=A_{i+1,i}=-1$
(we also introduce an insignificant modification $A_{-L,L}=A_{L,-L}=-1$ to have the 
cyclic structure of $A$), $\Delta \Phi ^{(y)}$ is the correction for membrane's points. 
Matrix $A$ has eigenvalues 
$$\lambda _h=4-2\cos \theta _h, \hskip 5mm \textrm{with} \hskip 2mm
\theta _h=2\pi h/(2L+1), \hskip 2mm  h\in \{ -L,...,L\}$$
and eigenvectors $V_h$ whose $l$-th component is simply $e^{il\theta _h}$.
We decompose $\Phi ^{(y)}(\theta )$ and $\Delta \Phi ^{(y)}(\theta )$ on the
base of $V_h$,
$$\Phi ^{(y)}(\theta )=\frac{1}{2L+1}\sum\limits _{h=-L}^L c_y(\theta ,\theta _h)V_h, 
\hskip 5mm  \Delta \Phi ^{(y)}(\theta )=\frac{1}{2L+1}\sum\limits _{h=-L}^L 
\Delta c_y(\theta ,\theta _h)V_h . $$
Their coefficients $c_y$ and $\Delta c_y$ obey the recurrence relations
\begin{equation}                                                               \label{eq:recurrence}
\lambda _hc_y(\theta ,\theta _h)=c_{y-1}(\theta ,\theta _h)+c_{y+1}(\theta ,\theta _h)+
\Delta c_y(\theta ,\theta _h) .
\end{equation}
%
So, the problem now is to find $c_y(\theta ,\theta _h)$.
It can be solved in two steps. First, using (\ref{eq:recurrence}) with
certain conditions, we can express $c_y$ in terms of $\Delta c_y$, $\lambda _h$ and 
$\varphi (\theta )$. Second, we have to find $\Delta c_y$. The first step was made in 
\cite{Grebenkov02}: 

-- taking sufficiently large $N_u$ and $N_l$, we impose two 
conditions to close the recurrence relations (\ref{eq:recurrence}),
\begin{equation}                                                               \label{eq:conditions}
c_{N_u+1}(\theta ,\theta _h)\approx \varphi (\theta )c_{N_u}(\theta ,\theta _h),  
\hskip 15mm  c_{-N_l}(\theta ,\theta _h)=0 .
\end{equation}
Note that the first condition is an approximate relation; 

-- we find the explicit solution of (\ref{eq:recurrence}), i.e., we express
$c_y$ in terms of $c_0$ and $\{\Delta c_{y'}\}$ (see Conclusions for more 
details);

-- we take the limit $N_u\to \infty $ and $N_l\to \infty $ to obtain
\begin{equation}                                                                      \label{eq:c_y}
c_y(\theta ,\theta _h)=\varphi ^{|y|}(\theta _h)c_0(\theta ,\theta _h)+
\sum\limits _{y'=-\N }^{N}\bigl[\gamma ^{(y)}_{y'}(\theta _h)+
\gamma ^{(-y)}_{-y'}(\theta _h)\bigr]\Delta c_{y'}(\theta ,\theta _h) ,
\end{equation}
where
\begin{equation}                                                                   \label{eq:gamma}
\gamma ^{(y)}_{y'}(\theta _h)=\sum\limits _{j=1}^{\min \{y,y'\} }
[\varphi (\theta _h)]^{2j-1+|y-y'|}   
\end{equation}
(here we use the convention that $\sum\nolimits _{j=a}^b$ is equal to $0$ if $b<a$, i.e.,
$\gamma ^{(y)}_{y'}$ is equal to $0$ if $y\leq 0$ or $y'\leq 0$). Note that
using only two simple identities,
\begin{equation}                                                                  \label{eq:two_iden}
\varphi ^2(\theta _h)-\lambda _h\varphi (\theta _h)+1=0,  \hskip 5mm
\gamma ^{(y+1)}_{y'}(\theta _h)-\lambda _h\gamma ^{(y)}_{y'}(\theta _h)+
\gamma ^{(y-1)}_{y'}(\theta _h)=-\delta _{y,y'}
\hskip 2mm  (y>0) ,
\end{equation}
one can verify directly, without intermediate steps of \cite{Grebenkov02}, 
that (\ref{eq:c_y}) is the solution of recurrence relations (\ref{eq:recurrence}).
Now we write $\phi _{x,y}(\theta )$ as
\begin{equation}                                                                   \label{eq:phi_xy}
\phi _{x,y}(\theta )=\sum\limits _{h=-L}^L \frac{e^{ix\theta _h}}{2L+1}
\biggl[\varphi ^{|y|}(\theta _h)c_0(\theta ,\theta _h)+\sum\limits _{y'=-\N }^{N}
\bigl(\gamma ^{(y)}_{y'}(\theta _h)+\gamma ^{(-y)}_{-y'}(\theta _h)\bigr)
\Delta c_{y'}(\theta ,\theta _h)\biggr] .
\end{equation}

So, the first step of our program is completed. Note that all the above
results are exactly the same as for a regular surface with compact support 
(by this reason we did not explain all details of calculation). Now we are 
ready to pass to the second step : to find $\Delta c_{y'}(\theta ,\theta _h)$ 
and $c_0(\theta ,\theta _h)$. These coefficients significantly depend on the 
membrane's geometry. Thus we cannot use the expressions of \cite{Grebenkov02}, 
and we should recalculate $\Delta c_{y'}(\theta ,\theta _h)$ 
and $c_0(\theta ,\theta _h)$.

\vskip 10mm
\section{ Two-dimensional finite membrane }                                         \label{sec:2D}
\subsection{ Coefficients $\Delta c_{y'}$. }                                        \label{sec:coeff}

To calculate $\Delta c_{y'}$ we should define explicitly what $\Delta \Phi ^{(y')}$ is.
We recall that these vectors were introduced in order to write correctly the 
expression (\ref{eq:relation}) for the membrane's points. In other words, 
relations (\ref{eq:relation}) are satisfied automatically for any external point, 
but they should be imposed artificially for any point of the membrane.

For membrane's point we have 
\begin{equation}                                                                     \label{eq:zeros}
\phi _{x,y}(\theta )=\cases{ e^{i\I (x,y)\theta }, \: \iff (x,y)\in \S ,\cr
                             \hskip 5mm  0, \hskip 7mm  \iff (x,y)\in \M\backslash \S, }
\end{equation}
because particles cannot penetrate in the depth of the membrane. The first line
of (\ref{eq:zeros}) represents the condition $P_{x_m,y_m}(n)=\delta _{m,n}$. Indeed,
if a particle starts from the boundary point $(x_m,y_m)$, it should be 
immediately absorbed by this point. In other words, the probability to be
absorbed by the $m$-th site is equal to $1$ while the other sites have no
chance. Expression (\ref{eq:zeros}) means that in $\Delta \Phi ^{(y')}$ there is only 
the contribution of boundary points and of membrane's points near boundary. 
A simple verification shows that 
\begin{equation}                                                                \label{eq:correction}
(\Delta \Phi ^{(y')})_{x'}=\chi _{\M}(x',y')\biggl[4e^{i\I (x',y')\theta }-
e^{i\I (x'+1,y')\theta }
\end{equation}
$$-e^{i\I (x'-1,y')\theta }-e^{i\I (x',y'+1)\theta }-e^{i\I (x',y'-1)\theta }-
\phi _{[\I (x',y')]}(\theta )\biggr] ,$$
where we enumerate near-boundary functions with the help of function $\I (x',y')$
(we use notation $\phi _{[m]}$ for the near-boundary function with index $m$). 
Again we insist to use our convention about $\I (x',y')$, i.e., in expression 
(\ref{eq:correction}) there are only the terms where corresponding point lies on $\S $. 

Using (\ref{eq:correction}), we determine $\Delta c_{y'}$ according to its definition,
$$\Delta c_{y'}(\theta ,\theta _h)=\sum\limits _{x'=-L}^L e^{-ix'\theta _h}(\Delta \Phi ^{(y')})_{x'} .$$

Usually there are a few nonzero components of $\Delta \Phi ^{(y')}$ for each $y'$.
Indeed, according to the formula (\ref{eq:correction}), $(\Delta \Phi ^{(y')})_{x'}$ is 
defined by the boundary points and by the membrane's points
near the boundary. But on the $n$-th level there are a few
such points, if $n\ne 0$. On the contrary, on the level zero 
there is an infinity of the boundary points due to the plane
``tails''. Thus, the vector $(\Delta \Phi ^{(-1)})$ has exceptional structure.
It contains the usual terms due to the non-trivial part of the membrane,
and the contribution of plane ``tails''.

\vskip 5mm
\subsection{ Boundary points' contribution }                                  \label{sec:boundary}

Let us calculate the contribution of the $m$-th boundary point to
$\phi _{x,y}(\theta )$ (see the second term in (\ref{eq:phi_xy})).
During these calculations we consider the case $y>0$, i.e., we write
just $\gamma ^{(y)}_{y'}$ omitting $\gamma ^{(-y)}_{-y'}$. 
The opposite case $y<0$ is obtained by reflection of all 
ordinates with respect to the horizontal axis, $y\to -y$,
$y'\to -y'$.

\begin{figure}             
\begin{center}
\includegraphics[width=75mm,height=5cm]{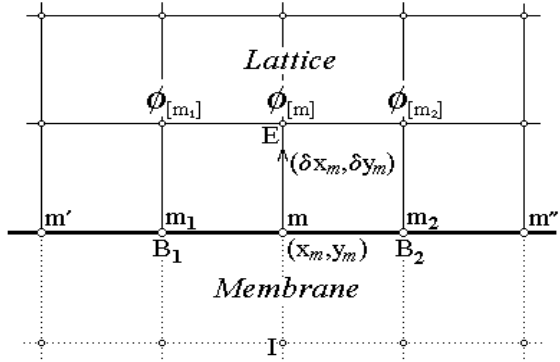}
\end{center}
\caption{ Boundary points' contribution. The boundary point $B=(x_m,y_m)$ 
with the outer normal $(\delta x_m,\delta y_m)=(0,1)$ is arounded by
two boundary points on the same level ($B_1$ and $B_2$), and by one 
external point $E$ and one internal points $I$.}                                \label{fig:boundary}
\end{figure}

Following the direction of the outer normal, there are four
possible positions of the boundary point. Let us consider 
the case when $(\delta x_m,\delta y_m)=(0,1)$ (see Fig.~\ref{fig:boundary}).
We suppose that the boundary point does not lie near the corners
(i.e., we consider {\it general position} of the boundary point). 
The points near corners will be considered separately in Appendices~\ref{sec:corner}. 
Note that these points give the additional corrections 
due to the lattice distortion caused by the procedure of corner's 
removing. Usually their influence is negligible, and we do not account 
them here (see Appendices for more details).

The boundary point $B=(x_m,y_m)$ has four neighbours: one external point
$E$, one internal point $I$, and two boundary points $B_1$ and $B_2$
with indices $m_1$ and $m_2$. Let us write accurately the corrections 
for all these points. Using (\ref{eq:correction}), we obtain
\begin{center}
\begin{tabular}{ c l l }
$B:$  &  $\hskip 4mm \Delta \phi _{x_m,y_m}(\theta )=4e^{im\theta }-\phi _{[m]}(\theta )
-e^{im_1\theta }-e^{im_2\theta }$,  & $e^{-ix_m\theta _h}\gamma ^{(y)}_{y_m}(\theta _h)$.  \\
$B_1:$        &  $\Delta \phi _{x_m-1,y_m}(\theta )=4e^{im_1\theta }-\phi _{[m_1]}(\theta )
-e^{im\theta }-e^{im'\theta }$,     & $e^{-i(x_m-1)\theta _h}\gamma ^{(y)}_{y_m}(\theta _h)$.  \\
$B_2:$        &  $\Delta \phi _{x_m+1,y_m}(\theta )=4e^{im_2\theta }-\phi _{[m_2]}(\theta )
-e^{im\theta }-e^{im''\theta }$,    & $e^{-i(x_m+1)\theta _h}\gamma ^{(y)}_{y_m}(\theta _h)$.    \\
$I:$          &  $\Delta \phi _{x_m,y_m-1}(\theta )=-e^{im\theta }$,   
                                 & $e^{-ix_m\theta _h}\gamma ^{(y)}_{y_m-1}(\theta _h)$.  \\
$E:$          &  $\Delta \phi _{x_m,y_m+1}(\theta )=0$,    
                                 & $e^{-ix_m\theta _h}\gamma ^{(y)}_{y_m+1}(\theta _h)$.\\
\end{tabular}
\end{center}
The third column represents the weight factor of corresponding point (see formula 
(\ref{eq:phi_xy})). Now we group together all terms containing factor $e^{im\theta }$,
and call such a group the boundary point's contribution $Z_m$ (we also add to this
group the corresponding near-boundary function $\phi _{[m]}$),
$$Z_m=e^{im\theta -ix_m\theta _h}\times \cases{ (4-2\cos \theta _h)\gamma ^{(y)}_{y_m}-
\gamma ^{(y)}_{y_m-1},  \hskip 13mm \iff  (\delta x_m,\delta y_m)=(0,1),    \cr
                                        (4-2\cos \theta _h)\gamma ^{(y)}_{y_m}-
\gamma ^{(y)}_{y_m+1},  \hskip 13mm \iff  (\delta x_m,\delta y_m)=(0,-1),     \cr
                                        (4-e^{-i \theta _h})\gamma ^{(y)}_{y_m}-
\gamma ^{(y)}_{y_m-1}-\gamma ^{(y)}_{y_m+1}, \hskip 2mm \iff  
 (\delta x_m,\delta y_m)=(-1,0),  \cr
                                        (4-e^{+i \theta _h})\gamma ^{(y)}_{y_m}-
\gamma ^{(y)}_{y_m-1}-\gamma ^{(y)}_{y_m+1}, \hskip 2mm \iff  
 (\delta x_m,\delta y_m)=(1,0)   }$$
$$-e^{-ix_m\theta _h}\gamma ^{(y)}_{y_m}\phi _{[m]}(\theta ) $$
(the first line corresponds to Fig.~\ref{fig:boundary}; three other cases can be 
obtained by analogy). We see that these contributions are similar. We complete each 
one by substracting and adding corresponding term, i.e.,
\begin{equation}                                                              \label{eq:supplem1}
Z_m=Z_m^{(complete)}+\cases{
e^{im\theta -ix_m\theta _h}\gamma ^{(y)}_{y_m+1}, \hskip 6mm \iff  
    (\delta x_m,\delta y_m)=(0,1),  \cr
e^{im\theta -ix_m\theta _h}\gamma ^{(y)}_{y_m-1}, \hskip 6mm \iff  
    (\delta x_m,\delta y_m)=(0,-1),  \cr
e^{im\theta -ix_m\theta _h}\gamma ^{(y)}_{y_m}e^{i\theta _h}, \hskip 4mm \iff  
    (\delta x_m,\delta y_m)=(-1,0), \cr
e^{im\theta -ix_m\theta _h}\gamma ^{(y)}_{y_m}e^{-i\theta _h}, \hskip 2mm \iff  
    (\delta x_m,\delta y_m)=(1,0),}
\end{equation}
where
$$Z_m^{(complete)}=e^{im\theta -ix_m\theta _h}\biggl[(4-2\cos \theta _h)\gamma ^{(y)}_{y_m}-
\gamma ^{(y)}_{y_m-1}-\gamma ^{(y)}_{y_m+1}\biggr]-
e^{-ix_m\theta _h}\gamma ^{(y)}_{y_m}\phi _{[m]}(\theta ). $$
Using the second identity of (\ref{eq:two_iden}), we write
$$Z_m^{(complete)}=e^{im\theta -ix_m\theta _h}\delta _{y,y_m}
-e^{-ix_m\theta _h}\gamma ^{(y)}_{y_m}\phi _{[m]}(\theta ).$$
With the help of the outer normal $(\delta x_m,\delta y_m)$, we can simplify (\ref{eq:supplem1}),
\begin{equation}                                                                     \label{eq:Z_m}
Z_m=e^{im\theta -ix_m\theta _h}\left(\delta _{y,y_m}+\gamma ^{(y)}_{y_m+\delta y_m}(\theta _h)
e^{-i\delta x_m\theta _h}\right)-e^{-ix_m\theta _h}\gamma ^{(y)}_{y_m}(\theta _h)\phi _{[m]}(\theta ).
\end{equation}

\vskip 5mm
\subsection { Distribution of hitting probabilities }                       \label{sec:distribution}

Now we can come back to the main problem. Expression (\ref{eq:phi_xy}) 
contains two terms which can be denoted as $\phi ^{(1)}_{x,y}(\theta )$ 
and $\phi ^{(2)}_{x,y}(\theta )$. We are going to simplify each of them independently.

Let us consider the first term,
$$\phi ^{(1)}_{x,y}(\theta )=\sum\limits _{h=-L}^L \frac{e^{ix\theta _h}}{2L+1}
\varphi ^{|y|}(\theta _h)c_0(\theta ,\theta _h) .$$
Following the definition of $c_y(\theta ,\theta _h)$, we have
$$c_0(\theta ,\theta _h)=(\Phi _0,V_h)=\sum\limits _{x'=-L}^Le^{-ix'\theta _h}
\phi _{x',0}(\theta ) .$$
The level zero can contain external and membrane's points.
Characteristic functions $\phi _{x,0}(\theta )$ for the 
membrane's points are given by (\ref{eq:zeros}). The external points on
the level zero were called ground points, and they are enumerated
by index $g=1..G$. Thus, 
$$c_0(\theta ,\theta _h)=\sum\limits _{x'=-L}^Le^{-ix'\theta _h}e^{i\I (x',0)\theta }
+\sum\limits _{g=1}^G e^{-i\bar{x}_g\theta _h}\phi _{\bar{x}_g,0}(\theta ).$$
Substituting this expression into $\phi ^{(1)}_{x,y}(\theta )$ and
taking the limit $L\to \infty $, we obtain
\begin{equation}                                                                     \label{eq:phi1}
\phi ^{(1)}_{x,y}(\theta )=\sum\limits _{x'=-\infty }^{\infty } 
e^{i\I (x',0)\theta }H^{y}_{x-x'}+\sum\limits _{g=1}^G H^{y}_{x-\bar{x}_g}
\phi _{\bar{x}_g,0}(\theta ) .
\end{equation}
This is the contribution of the level zero.

\vskip 1mm 

The second term,
$$\phi ^{(2)}_{x,y}(\theta )=\sum\limits _{h=-L}^L \frac{e^{ix\theta _h}}{2L+1}
\sum\limits _{y'=-\N }^{N}\sum\limits _{x'=-L}^L
\gamma ^{(y)}_{y'}(\theta _h) e^{-ix'\theta _h} (\Delta \Phi ^{(y')})_{x'} ,$$
represents other levels except level zero. Summation over $x'$ and $y'$ can be
replaced by the sum of the boundary points' contributions,
$$\phi ^{(2)}_{x,y}(\theta )=\sum\limits _{h=-L}^L \frac{e^{ix\theta _h}}{2L+1}
\sum\limits _{m=1}^MZ_m .$$
Using (\ref{eq:Z_m}), we write the limit $L\to \infty $ for the case $y\geq 0$,
$$\phi ^{(2)}_{x,y}(\theta )=\int\limits _{-\pi }^{\pi }\frac{d\theta '}{2\pi }
\sum\limits _{m=1}^Me^{i(x-x_m)\theta '}\biggl[e^{im\theta }\delta _{y,y_m}+
e^{im\theta }\gamma ^{(y)}_{y_m+\delta y_m}(\theta ')e^{-i\delta x_m\theta '}-
\gamma ^{(y)}_{y_m}(\theta ')\phi _{[m]}(\theta )\biggr] .$$
Changing the order of summation and integration, we obtain
\begin{equation}                                                                     \label{eq:phi2}
\phi ^{(2)}_{x,y}(\theta )=\sum\limits _{m=1}^M \biggl[e^{im\theta }
\delta _{x,x_m}\delta _{y,y_m}+e^{im\theta }D^{y,y_m+\delta y_m}_{x-x_m-\delta x_m}-
D^{y,y_m}_{x-x_m}\phi _{[m]}(\theta )\biggr]
\end{equation}
where we introduced coefficients
$$D^{y,y'}_x=\int\limits _{-\pi }^{\pi }\frac{d\theta }{2\pi }e^{ix\theta }
\bigl(\gamma ^{(y)}_{y'}(\theta )+\gamma ^{(-y)}_{-y'}(\theta )\bigr) .$$
Note that function $\gamma ^{(-y)}_{-y'}(\theta )$ is added into the
definition of $D^{y,y'}_x$ in order to generalize it for the case $y<0$.
It means that formula (\ref{eq:phi2}) is already valid for the both cases $y\geq 0$
and $y\leq 0$. Using definitions (\ref{eq:H}) and (\ref{eq:gamma}), we rewrite 
this expression in terms of $H^y_x$, 
\begin{equation}                                                                         \label{eq:D}
D^{y,y'}_x=\cases{ \sum\limits _{j=1}^{\min \{|y|,|y'|\}}H^{2j-1+|y-y'|}_x , 
\hskip 5mm  \iff y\cdot y'>0 ,\cr
    \hskip 16mm              0, \hskip 23mm  \iff  y\cdot y'\leq 0 . }
\end{equation}

We bring expressions (\ref{eq:phi1}) and (\ref{eq:phi2}) together to obtain
\begin{equation}                                                                       \label{eq:phi}
\phi _{x,y}(\theta )=\tilde{\phi }_{x,y}(\theta )-\sum\limits _{m=1}^M
D^{y,y_m}_{x-x_m}\phi _{[m]}(\theta )+\sum\limits _{g=1}^G
H^{y}_{x-\bar{x}_g}\phi _{\bar{x}_g,0}(\theta ) ,
\end{equation}
where
\begin{equation}                                                                 \label{eq:phi_tilde}
\tilde{\phi }_{x,y}(\theta )=\sum\limits _{m=-\infty }^{\infty }
e^{im\theta }\left(\delta _{x,x_m}\delta _{y,y_m}+
D^{y,y_m+\delta y_m}_{x-x_m-\delta x_m}\right) .
\end{equation}

Here we should explain how (\ref{eq:phi_tilde}) is obtained. This sum over $m$
can be separated in two parts: the sum over $m\in [1,M]$ and the sum 
over $m\notin [1,M]$. The first part contained in (\ref{eq:phi2}) corresponds
to the boundary points lying on all levels except level zero. The second part
gives contributions of boundary points on the level zero, including
the plane ``tails''. In the case $y>0$, expression (\ref{eq:phi1}) contains the sum
$$\sum\limits _{x'=-\infty }^{\infty }e^{i\I (x',0)\theta }H^{y}_{x-x'} ,$$
which can be represented as
$$\sum\limits _{m\notin [1,M]}e^{im\theta }H^y_{x-x_m}.$$
For such points $y_m=0$, $\delta x_m=0$. Supposing $\delta y_m=1$
and using relation $D^{y,1}_x=H^y_x$ (for $y>0$), we immediatly obtain
$$\sum\limits _{m\notin [1,M]}e^{im\theta }D^{y,y_m+\delta y_m}_{x-x_m-\delta x_m} ,$$
that proves representation (\ref{eq:phi_tilde}) for $\tilde{\phi }_{x,y}(\theta )$.
Now, if $y<0$, we {\it should not} write the contribution of plane
``tails'' because it is compensated by $\Delta c_{-1}$ (see the end of Section~\ref{sec:coeff}).
In other words, the plane ``tails'' have no {\it direct} influence on
the points in the lower half plane\footnote{ See Conclusions 
for more detailed discussion on this topic.}. The same concerns
the case with $y>0$ and $\delta y_m<0$. In the case $y=0$ coefficient
$H^y_x$ becomes $\delta $-symbol which appears in expression (\ref{eq:phi_tilde})
explicitly. We conclude that expressions (\ref{eq:phi}) and (\ref{eq:phi_tilde}) 
are valid for any external point $(x,y)$, in particular, with $y<0$ 
(if such points exist). 

\vskip 1mm

Applying the inverse Fourier transform to (\ref{eq:phi}), we obtain the
distribution of hitting probabilities,
\begin{equation}                                                                        \label{eq:P}
P_{x,y}(n)=\delta _{x,x_n}\delta _{y,y_n}+D^{y,y_n+\delta y_n}_{x-x_n-\delta x_n}
-\sum\limits _{m=1}^M D^{y,y_m}_{x-x_m}P_{[m]}(n)+\sum\limits _{g=1}^G
H^{y}_{x-\bar{x}_g}P_{\bar{x}_g,0}(n) .
\end{equation}

\vskip 5mm
\subsection{ Equations for near-boundary and ground functions }               \label{sec:equations}

To complete our calculations, we should find the near-boundary functions
$P_{[m]}(n)$ and ground functions $P_{\bar{x}_g,0}(n)$ entering in expression 
(\ref{eq:P}). After that, one can use this expression for any $x$, $y$ and $n$.
We take  
$$x=x_k+\delta x_k, \hskip 5mm y=y_k+\delta y_k  \hskip 5mm \textrm{for} \hskip 5mm 
k\in [1,M]$$
to obtain $M$ linear equations for the near-boundary functions $P_{[k]}(n)$,
$k\in [1,M]$
\begin{equation}                                                                  \label{eq:supplem2}
P_{[k]}(n)=D^{y_k+\delta y_k,y_n+\delta y_n}_{x_k+\delta x_k-x_n-\delta x_n}-
\sum\limits _{m=1}^M D^{y_k+\delta y_k,y_m}_{x_k+\delta x_k-x_m}P_{[m]}(n)
+\sum\limits _{g=1}^G H^{y_k+\delta y_k}_{x_k+\delta x_k-\bar{x}_g}P_{\bar{x}_g,0}(n) .
\end{equation}
Let us introduce matrices,
$$(D_{NN})_{k,m}=D^{y_k+\delta y_k,y_m}_{x_k+\delta x_k-x_m},  
\hskip 5mm   k\in [1,M],\: m\in [1,M]$$  
$$(D_{NG})_{k,m}=H^{y_k+\delta y_k}_{x_k+\delta x_k-\bar{x}_m},  
\hskip 5mm  k\in [1,M],\: m\in [1,G].$$
If there is no ground functions, i.e., the non-trivial part of 
the membrane completely lies in the upper half plane, we take $D_{NG}=0$. 

Now we can rewrite (\ref{eq:supplem2}) as
\begin{equation}                                                                  \label{eq:supplem3}
(I+D_{NN})\left(\begin{array}{ c } P_{[1]}(n) \\ ... \\ P_{[M]}(n) \\ \end{array}  
\right) -D_{NG}\left(\begin{array}{ c } P_{\bar{x}_1,0}(n) \\ ... \\ P_{\bar{x}_G,0}(n) \\ 
\end{array}  \right) = \left(\begin{array}{ c } 
D^{y_1+\delta y_1,y_n+\delta y_n}_{x_1+\delta x_1-x_n-\delta x_n} \\ ... \\ 
D^{y_M+\delta y_M,y_n+\delta y_n}_{x_M+\delta x_M-x_n-\delta x_n} \\  \end{array}\right). 
\end{equation}

If ground points exist, we take the discrete Laplacian equations (\ref{eq:Laplace_discrete})
for these points in order to obtain conditions for ground functions,
$$4P_{\bar{x}_g,0}(n)-P_{\bar{x}_g-1,0}(n)-P_{\bar{x}_g+1,0}(n)-P_{\bar{x}_g,1}(n)-
P_{\bar{x}_g,-1}(n)=0, \hskip 3mm  g\in [1..G].$$
Substituting $P_{\bar{x}_g,1}(n)$ and $P_{\bar{x}_g,-1}(n)$ from (\ref{eq:P})
into these conditions, we obtain
$$4P_{\bar{x}_g,0}(n)-P_{\bar{x}_g+1,0}(n)-P_{\bar{x}_g-1,0}(n)=
D^{1,y_n+\delta y_n}_{\bar{x}_g-x_n-\delta x_n}+D^{-1,y_n+\delta y_n}_{\bar{x}_g-x_n-\delta x_n}$$
$$-\sum\limits _{m=1}^M \bigl(D^{1,y_m}_{\bar{x}_g-x_m}+D^{1,-y_m}_{\bar{x}_g-x_m}\bigr)
P_{[m]}(n)+\sum\limits _{g'=1}^G 2H^{1}_{\bar{x}_g-\bar{x}_{g'}}P_{\bar{x}_{g'},0}(n).$$
In matrix form, 
\begin{equation}                                                                  \label{eq:supplem4}
D_{GN} \left(\begin{array}{ c } P_{[1]}(n) \\ ... \\ 
P_{[M]}(n)   \\ \end{array}  \right) +
D_{GG}\left(\begin{array}{c}  P_{\bar{x}_1,0}(n) \\ ... \\ 
P_{\bar{x}_G,0}(n)  \\  \end{array} \right) = \frac12
\left(\begin{array}{c} H^{y_n+\delta y_n}_{\bar{x}_1-x_n-\delta x_n} \\ ... \\
H^{y_n+\delta y_n}_{\bar{x}_G-x_n-\delta x_n} \\ \end{array} \right) ,
\end{equation}
where
$$(D_{GG})_{k,m}=-H^{1}_{\bar{x}_k-\bar{x}_m}+
\cases{ \hskip 3mm  2, \hskip 5mm  \iff  k=m, \cr
       -0.5, \hskip 3mm \iff |\bar{x}_k-\bar{x}_m|=1,  \hskip 5mm  k\in [1,G], \: m\in [1,G], \cr
       \hskip 3mm  0, \hskip 6mm  \textrm{otherwise} }$$
$$(D_{GN})_{k,m}=\frac12H^{y_m}_{\bar{x}_k-x_m}, \hskip 50mm  
k\in [1,G], \: m\in [1,M] .$$

We have two matrix equations (\ref{eq:supplem3}) and (\ref{eq:supplem4}) which
allow to find near-boundary and ground functions,
\begin{equation}                                                                   \label{eq:matrix1}
\left(\begin{array}{ c } P_{[1]}(n) \\ ... \\ P_{[M]}(n) \\ \end{array}  \right)=
(I+D)^{-1}\left(\begin{array}{ c } P^*_{[1]}(n) \\ ... \\ 
P^*_{[M]}(n) \\ \end{array}\right)
\end{equation}
where
$$D=D_{NN}+D_{NG}D_{GG}^{-1}D_{GN} \hskip 5mm   \textrm{and}$$
\begin{equation}                                                                        \label{eq:P*}
\left(\begin{array}{ c } P^*_{[1]}(n) \\ ... \\ P^*_{[M]}(n) \\ \end{array}\right)
=\left(\begin{array}{ c } D^{y_1+\delta y_1,y_n+\delta y_n}_{x_1+\delta x_1-x_n-\delta x_n} \\ ... \\ 
D^{y_M+\delta y_M,y_n+\delta y_n}_{x_M+\delta x_M-x_n-\delta x_n}  \\ \end{array}\right)+\frac12D_{NG}D_{GG}^{-1}
\left(\begin{array}{c} H^{y_n+\delta y_n}_{\bar{x}_1-x_n-\delta x_n} \\ ... \\
H^{y_n+\delta y_n}_{\bar{x}_G-x_n-\delta x_n} \\ \end{array} \right) .
\end{equation}
Using (\ref{eq:supplem4}) again, we obtain ground functions,
\begin{equation}                                                                    \label{eq:matrix2}
\left(\begin{array}{ c } P_{\bar{x}_1,0}(n) \\ ... \\  P_{\bar{x}_G,0}(n) \\
\end{array}\right)=\frac12D_{GG}^{-1}\left(\begin{array}{c}  
H^{y_n+\delta y_n}_{\bar{x}_1-x_n-\delta x_n} \\ ... \\
H^{y_n+\delta y_n}_{\bar{x}_G-x_n-\delta x_n}  \\ \end{array} \right)-D_{GG}^{-1}D_{GN}
(I+D)^{-1}\left(\begin{array}{ c } P^*_{[1]}(n) \\ ... \\ P^*_{[M]}(n)
\end{array}\right) .
\end{equation}

Expression (\ref{eq:P}) is our main result. What have we done?
Using the characteristic functions technique, we express
the hitting probability $P_{x,y}(n)$ (for any $x$, $y$ and $n$) in
terms of the explicit coefficients $D^{y,y'}_x$ and {\em finite}
number of coefficients $P_{[m]}(n)$ and $P_{\bar{x}_g,0}(n)$
which can be calculated with the help of (\ref{eq:matrix1}), (\ref{eq:P*})
and (\ref{eq:matrix2}).

\vskip 1mm

  Note that $P_{x,y}(n)$ depends on $n$ only through coordinates $x_n$ and $y_n$ 
of the $n$-th boundary point and through outer normal $(\delta x_n,\delta y_n)$ at 
this point. It means that obtained distribution of hitting probabilities 
{\it does not depend} on a choice of parametrization of the membrane.
In other words, we can use any parametrization of the membrane's 
boundary. This remark will be used in the next section.

\vskip 10mm
\section{  Multi-dimensional membranes }                                         \label{sec:multi}

We confined ourselves to the two-dimensional case in order to 
clarify the treatment. Now we are going to generalize the 
previous results for multi-dimensional membranes.
  
Consider $d$-dimensional hypercubic lattice. As above, we can
define a finite membrane with the help of similar conditions as in 
Section~\ref{sec:definitions}. In particular, we suppose that membranes have a 
compact support (the compactness condition), i.e., there exists a $(d-1)$-dimensional
hyperplane (level zero) such that the membrane's boundary is just
its ``finite'' perturbation. We choose such coordinates that 
this hyperplane is defined by equation $x^{(d)}=0$, its points
are enumerated by multi-index $\bx =(x^{(1)},...,x^{(d-1)})$. 
Coordinate $x^{(d)}$ of the orthogonal direction is denoted $y$.
Let us denote $e_k$ the unit vectors in $k$-th direction, $k\in \{1,...,d-1\}$.
  
As usual, we enumerate the boundary points by index $n$,
$(\bx _n,y_n)$. The outer normal on $n$-th point is
$(\delta \bx _n,\delta y_n)$.
We can introduce the $n$-th level as a set of points
lying on the hyperplane with $y=n$, i.e., $\{ (\bx ,y)\: :\: y=n, \: \forall \: \bx \}$. 

The discrete Laplacian equation in $d$-dimensional case is
$$P_{\bx ,y}(n)=\frac{1}{2d}\sum\limits _{i=1}^{d-1}\biggl(P_{\bx +e_i,y}(n)
+P_{\bx -e_i,y}(n)\biggr)+\frac{1}{2d}\biggl(P_{\bx ,y+1}(n)+P_{\bx ,y-1}(n)\biggr) ,$$ 
whence one can obtain the distribution of hitting probabilities 
for a planar multi-dimensional case,
\begin{equation}                                                                       \label{eq:H_d}
H^{y}_{\bx }=\int\limits _{-\pi }^{\pi }...\int\limits _{-\pi }^{\pi }
\frac{d\theta _1...d\theta _{d-1}}{(2\pi )^{d-1}}e^{i(\bx \cdot \bt )}\varphi ^{|y|}(\theta _i)
\end{equation}
where $(\bx \cdot \bt )$ is a scalar product of two $(d-1)$-dimensional
vectors $\bx $ and $\bt $; function $\varphi (\theta _i)$ depends on 
$\theta _1$,..., $\theta _{d-1}$,
$$\varphi (\theta _1,...\theta _{d-1})=\biggl(d-\sum\limits _{i=1}^{d-1}\cos \theta _i\biggr)
-\sqrt{\biggl(d-\sum\limits _{i=1}^{d-1}\cos \theta _i\biggr)^2-1} .$$

The crucial idea of our generalization is that all manipulations 
along the orthogonal direction remain valid. In particular,
the expressions (\ref{eq:c_y}) and (\ref{eq:gamma}) are true, if one takes 
$\varphi (\theta _1,...,\theta _{d-1})$ instead of $\varphi (\theta )$. 
It means that we can impose the same conditions (\ref{eq:conditions})
basing again on the technique proposed in \cite{Grebenkov02}. 
Now one can repeat the previous calculations of Section~\ref{sec:2D}:

-- to calculate coefficients $\Delta c_y$;

-- to obtain contributions of boundary points;

-- to impose conditions for near-boundary and ground functions.

Having made these technical steps, we understand that the only coefficients
$H^y_x$ (and $D^{y,y'}_x$ as consequence) are different for $d>2$, but the
structure of solution is exactly the same. It means that we leave
expression (\ref{eq:P}) without changes for $d$-dimensional case,
\begin{equation}                                                                       \label{eq:P_d}
P_{\bx ,y}(n)=\delta _{\bx ,\bx _n}\delta _{y,y_n}+
D^{y,y_n+\delta y_n}_{\bx -\bx _n-\delta \bx _n}-
\sum\limits _{m=1}^M D^{y,y_m}_{\bx -\bx _m}P_{[m]}(n)
+\sum\limits _{g=1}^G H^{y}_{\bx -\bar{\bx }_g}P_{\bar{\bx }_g,0}(n).
\end{equation}
Here coefficients $D^{y,y_n}_{\bx-\bx _n}$ are still defined with the help
of (\ref{eq:D}) but for coefficients $H^y_{\bx -\bx '}$ we should use 
(\ref{eq:H_d}) instead of (\ref{eq:H}).
 
As above, we can obtain a set of linear equations for near-boundary and 
ground functions. We take
$$P_{[k]}(n)=D^{y_k+\delta y_k,y_n+\delta y_n}_{\bx _k+\delta \bx _k-\bx _n-\delta \bx _n}-
\sum\limits _{m=1}^M D^{y_k+\delta y_k,y_m}_{\bx _k+\delta \bx _k-\bx _m}P_{[m]}(n)
+\sum\limits _{g=1}^G H^{y_k+\delta y_k}_{\bx _k+\delta \bx _k-\bar{\bx }_g}P_{\bar{\bx }_g,0}(n)$$
for near-boundary functions, and
$$P_{\bar{\bx }_g,0}(n)=\frac{1}{2d}\sum\limits _{i=1}^{d-1}\bigl(P_{\bar{\bx }_g+e_i,0}(n)+
P_{\bar{\bx }_g-e_i,0}(n)\bigr)
+\frac{1}{2d}\bigl(P_{\bar{\bx }_g,1}(n)+P_{\bar{\bx }_g,-1}(n)\bigr) $$
for ground functions. Using the same representation as in Section~\ref{sec:equations}, 
we obtain the near-boundary and ground functions according to expressions (\ref{eq:matrix1}) 
and (\ref{eq:matrix2}). The only change that we should make is concerned with matrix $D_{GG}$. 
Indeed, we generalize this matrix as
$$(D_{GG})_{k,m}=-H^{1}_{\bar{\bx }_k-\bar{\bx }_m}+
\cases{ \hskip 2mm  d, \hskip 6mm  \iff  k=m, \cr
       -0.5, \hskip 3mm \iff |\bar{\bx }_k-\bar{\bx }_m|=1, \cr
       \hskip 2mm  0, \hskip 7mm  \textrm{otherwise.} }$$

It is important to stress that we used the remark at the end of Section~\ref{sec:2D}: 
the distribution of hitting probabilities should not depend on a
choice of parametrization. It allows to avoid a complex multi-dimensional 
parametrization of the membrane. To be more rigorous, one should introduce
such a parametrization, recalculate distribution $P_{\bx ,y}$ again, and 
re-enumerate the boundary points by the single index in order to obtain 
(\ref{eq:P_d}). We omit these technical details. 

Note that expression (\ref{eq:P_d}) has the same structure for any 
dimension of the lattice, the only distinction is contained in
coefficients $H^{y}_{\bx }$ which are individual for each $d$. 
We give a useful asymptotics for coefficients $H^y_{\bx }$,
\begin{equation}                                                                   \label{eq:asympt}
H^{y}_{\bx }\approx \frac{\Gamma (d/2)}{\pi ^{d/2}}\frac{|y|}{(\bx ^2+y^2)^{d/2}} 
\end{equation}
that is the multi-dimensional generalization of the well-known 
Cauchy distribution of hitting probabilities of the Brownian motion.

\vskip 10mm
\section{ Brownian self-transport operator }                                        \label{sec:BSO}
\vskip 2mm

As it was mentioned before, we are interested
in the problems of Laplacian transfer, and consequently, in
the Brownian self-transport operator $Q$. We recall that $Q_{k,n}$ 
is the probability that a random walker contacts at the first hit 
the $n$-th site of the membrane if started from the $k$-th site,
without touching the other sites. Obviously, any path starting
from $k$-th site of membrane must pass through the
corresponding near-boundary point. Provided that we removed
all the singular (corner) points, any boundary point
has exactly one near-boundary point. It means that
$$Q_{k,n}=P_{[k]}(n) ,$$
i.e., the matrix $Q$ is composed of the near-boundary functions
$P_{[k]}(n)$. If one denotes
$$Q^*_{k,n}=P^*_{[k]}(n)  ,$$
expression (\ref{eq:matrix1}) becomes
\begin{equation}                                                                        \label{eq:Q}
Q=(I+D)^{-1}Q^* .
\end{equation}
Recall that the system of linear equations (\ref{eq:supplem2}) was obtained
only for non-trivial part of the membrane, in particular, matrix $Q$ 
has $M\times M$ elements. However, we are usually interested in the
whole Brownian self-transport operator, including the plane ``tails''.

To construct the Brownian self-transport operator for all sites
of our interest, we can slightly modify definitions of matrices
$D_{NN}$, $D_{NG}$ and $D_{GN}$. Let $\A $ be the set of indices
of sites for which one wants define $Q$. For example, if one
needs to know $Q$ for all sites of the membrane, it is sufficient
to take $\A =\Z $. For the numerical simulations one takes 
a finite number of sites. In any case, $\A $ must 
contain all non-trivial sites, i.e., $[1..M]\subset \A $. 
We remark that summation over $m$ in (\ref{eq:P}) can be elarged
to all integer numbers from $-\infty $ to $\infty $. Indeed,
for $m\notin [1,M]$ we have $y_m=0$, and $D^{y,y_m}_{x-x_m}=0$
according to (\ref{eq:D}). In particular, this summation covers
all possible sites which can be of our interest.
Now we rewrite the definitions of matrices $D_{NN}$, $D_{NG}$ and $D_{GN}$,
\begin{eqnarray*}
(D_{NN})_{k,m} &=D^{y_k+\delta y_k,y_m}_{x_k+\delta x_k-x_m},  &
\hskip 5mm   k\in \A,\hskip 3mm m\in \A    \\
(D_{NG})_{k,m} &=H^{y_k+\delta y_k}_{x_k+\delta x_k-\bar{x}_m},  &
\hskip 5mm  k\in \A,\: m\in [1,G]. \\
(D_{GN})_{k,m} &=\frac12H^{y_m}_{\bar{x}_k-x_m},  & \hskip 5mm  
k\in [1,G], \: m\in \A . \\
\end{eqnarray*}
Matrix $D_{GG}$ remains unchanged, and
$$D=D_{NN}+D_{NG}D_{GG}^{-1}D_{GN} .$$
Evidently, matrix $Q^*$ (i.e., vectors $P^*_{[m]}$) should be
also recalculated for all sites of the interest.

These modifications were made to obtain the Brownian self-transport 
operator in the matrix form (\ref{eq:Q}). We stress that after
finding the near-boundary and ground functions, one can use
(\ref{eq:P}) directly to ``compose'' the matrix $Q$,  
$$Q_{k,n}=D^{y_k+\delta y_k,y_n+\delta y_n}_{x_k+\delta x_k-x_n-\delta x_k}-
\sum\limits _{m=1}^M D^{y_k+\delta y_k,y_m}_{x_k+\delta x_k-x_m}P_{[m]}(n)
+\sum\limits _{g=1}^G H^{y_k+\delta y_k}_{x_k-\bar{x}_g}P_{\bar{x}_g,0}(n) .$$

\vskip 1mm

We can conclude that the Brownian self-transport operator is explicitly 
constructed for any reasonable multi-dimensional {\em finite} membrane. 
Now it can be useful to study the general properties of $Q^*$ 
and the spectral properties of $D$, and their dependencies on the 
geometry of the membrane. This work is not yet finished.

\vskip 10mm
\section{ Generalizations }                                                \label{sec:generalizations}

As we explained in Section~\ref{sec:multi}, expression (\ref{eq:P}) can be
generalized for the multi-dimensional case. One can go further
by regarding more complex problems. In this section
we consider the influence of different barriers, and the time-dependent
distributions of hitting probabilities.

It is important to stress that we have obtained the Brownian
self-transport operator for a membrane with infinite ``tails''
placed in the infinite space. However, usual physical membranes 
have a finite size, and they are placed in a closed volume 
(like a box or a cell). In order to apply our results to real membranes,
these features should be taken into account, for example, by
introducing absorbing or reflecting barriers.
First, we consider the horizontal barriers, and then
the vertical barriers.

\vskip 2mm
\subsection{ Horizontal barriers }                                           \label{sec:horizontal}

  Suppose that there exists a horizontal barrier at the 
level $\NN +1>N$ which absorbs or reflects all particles.
In this case, it is not difficult to recalculate 
the distribution of hitting probabilities (\ref{eq:P}).
Let us return to the recurrence relations (\ref{eq:recurrence}).
In order to close them, one can use a condition like
$$c_{\NN +1}(\theta ,\theta _h)=\eta (\theta ,\theta _h)
c_{\NN }(\theta ,\theta _h)$$
with a certain coefficient $\eta $ which can depend on 
$\theta $ and $\theta _h$. Now one can express $c_y$
in terms of $c_0$ and $\{ \Delta c_{y'}\}$ (and of 
some explicit functions like $\eta $ and $\varphi $). 
This expression was found in \cite{Grebenkov02} ,
$$c_y=f^{(\NN)}_y\biggl(c_0+\sum\limits _{y'=1}^{y}\alpha _{y'}\Delta c_{y'}\biggr)
+\alpha _y\sum\limits _{y'=y+1}^{\NN}f^{(\NN)}_{y'}\Delta c_{y'},$$
where 
$$f^{(\NN)}_y=\frac{[1-\eta \varphi (\theta _h)]-\varphi ^{2(\NN -y)+1}(\theta _h)
[\varphi (\theta _h)-\eta ]}{[1-\eta \varphi (\theta _h)]-\varphi ^{2\NN +1}(\theta _h)
[\varphi (\theta _h)-\eta ]}\varphi ^y(\theta _h) $$  
and
$$\alpha _y(\theta _h)=\sum\limits _{j=0}^{y-1}[\varphi (\theta _h)]^{2j+1-y}.$$
When level $\NN $ of the barrier goes to infinity, function
$f^{(\NN)}_y$ tends to $f^{(\infty )}_y=\varphi ^y$, independently of the
value of $\eta $. It simply means that one can choose any condition 
at infinity. 

The absorbing barrier corresponds to $\eta =0$ since the
particles cannot reach the membrane from any point $(x,\NN+1)$; thus
$$f^{(\NN, abs)}_y=\frac{1-\varphi ^{2\NN+2-2y}}{1-\varphi ^{2\NN+2}}\varphi ^y .$$
The reflecting barrier corresponds to $\eta =1$ since the
particles should have the same probability to reach the membrane
from points $(x,\NN+1)$ and $(x,\NN )$; thus
$$f^{(\NN, ref)}_y=\frac{1+\varphi ^{2\NN+1-2y}}{1+\varphi ^{2\NN+1}}\varphi ^y .$$
All other calculations are left unchanged, but we should replace
$\varphi ^y$ by $f^{(\NN, abs)}_y$ or $f^{(\NN, ref)}_y$. It means
that one can use expression (\ref{eq:P}) where coefficients
$H^y_x$ are replaced by modified ones,
$$\tilde{H}^y_x=\cases{   \hskip 15mm  0  \hskip 15mm  \iff y>\NN, \cr 
\int\limits _{-\pi }^{\pi }\frac{d\theta '}{2\pi }e^{ix\theta '}
f^{(\NN )}_y(\theta ')  \hskip 3mm \iff 1\leq y\leq \NN , \cr
 \hskip 13mm   H^{-y}_x  \hskip 11.5mm  \iff y<0   .}$$
The last line represents the fact that the barrier in the upper half plane
has no influence on the lower half plane.

Using this definition of coefficients $\tilde{H}^y_x$, one can
calculate the distribution of hitting probabilities for a 
general membrane with an additional barrier absorbing or
reflecting particles. Note that these modifications are valid
for any dimension.

\vskip 5mm
\subsection{ Vertical barriers }                                           \label{sec:vertical}

Now we can limit the particle's movement along the horizontal axis
putting an absorbing or reflecting vertical barrier. Let us put two 
barriers at $x=-L-1$ and $x=L+1$ so that the non-trivial part 
of the membrane lies between these vertical lines. In Section~\ref{sec:technique},
we represented the Laplacian equations (\ref{eq:relation}) in 
matrix form using matrix $A$,
$$A_{i,i}=4, \hskip 3mm  A_{i,i+1}=A_{i+1,i}=-1,  \hskip 3mm A_{-L,L}=A_{L,-L}=-1 ,$$
and all other elements are $0$. The last equality was imposed artificially
to obtain the cyclic structure of matrix $A$ that was convenient to 
have explicit eigenvectors $V_h$. After all calculations,
we took the limit $L\to \infty $, therefore this modification
vanished. Now we are working with {\it finite} $L$, consequently,
we should be more accurate with these conditions. Note that
relation $A_{-L,L}=A_{L,-L}=-1$ corresponds to {\it periodic or cyclic
boundary condition}, i.e., we maintain $\phi _{L+1,y}=\phi _{-L,y}$ and
$\phi _{-L-1,y}=\phi _{L,y}$. In other words, we suppose that
each pair of points $(L+1,y)$ and $(-L,y)$ is identified or ``glued together''
(the same for points $(-L-1,y)$ and $(L,y)$).

In order to introduce vertical barriers, one should change the
boundary conditions. For two absorbing barriers one should simply
remove the artificial condition, i.e., $A_{-L,L}=A_{L,-L}=0$
given that $\phi _{L+1,y}=\phi _{-L-1,y}=0$ (particles are absorbed
by barriers). For two reflecting barriers one should also
introduce the following modifications, $A_{-L,-L}=A_{L,L}=3$
given that $\phi _{L+1,y}=\phi _{L,y}$ and $\phi _{-L-1,y}=\phi _{-L,y}$.
Obviously, one can also consider the mixed case with one absorbing barrier
and one reflecting barrier. Eigenvalues and eigenvectors of these modified 
matrices have no explicit form. It is more convenient to use the
initial matrix $A$ (cyclic boundary conditions) with additional
corrections, i.e.,
$$A\Phi ^{(y)}=\Phi ^{(y-1)}+\Phi ^{(y+1)}+\Delta \Phi ^{(y)}+\Delta \tilde{\Phi }^{(y)},$$
where we introduced a new correction vector $\Delta \tilde{\Phi }^{(y)}$ whose components
are equal to $0$ except  
$$\begin{array}{ l l c }
(\Delta \tilde{\Phi }^{(y)})_{-L}=-\phi _{L,y} &
(\Delta \tilde{\Phi }^{(y)})_{L}=-\phi _{-L,y} & \textrm{(absorption)} \\ 
(\Delta \tilde{\Phi }^{(y)})_{-L}=-\phi _{L,y}+\phi _{-L,y} & 
(\Delta \tilde{\Phi }^{(y)})_{L}=-\phi _{-L,y}+\phi _{L,y}  &  
\textrm{(reflection)}. \\ \end{array}$$
Now we can repeat the previous calculations to obtain a new distribution
of hitting probabilities $P_{x,y}(n)$ instead of (\ref{eq:P}),
\begin{equation}                                                                  \label{eq:supplem5}
P_{x,y}(n)=\tilde{D}^{y,y_n+\delta y_n}_{x-x_n-\delta x_n}-
\sum\limits _{m=1}^M \tilde{D}^{y,y_m}_{x-x_m}P_{[m]}(n)+\sum\limits _{g=1}^G 
\tilde{H}^y_{x-\bar{x}_g} P_{\bar{x}_g,0}(n)-
\end{equation}
$$\cases{ \hskip 3mm  \sum\limits _{y'=1}^{\infty }\bigl[\tilde{D}^{y,y'}_{x-L}P_{L,y'}(n)+
\tilde{D}^{y,y'}_{x+L}P_{-L,y'}(n) \bigr],  \hskip 21mm   \textrm{(absorption)} ,\cr
\hskip 3mm  \sum\limits _{y'=1}^{\infty }(\tilde{D}^{y,y'}_{x-L}-\tilde{D}^{y,y'}_{x+L})
\bigl[P_{L,y'}(n)-P_{-L,y'}(n) \bigr], \hskip 16mm  \textrm{(reflection)} ,} $$
where $\tilde{H}^y_x$ should be calculated as {\it finite} sum over $h$
from $-L$ to $L$,
$$\tilde{H}^y_x=\sum\limits _{h=-L}^L \frac{e^{ix\theta _h}}{2L+1}\varphi ^{|y|}(\theta _h),
\hskip 10mm  \theta _h=\frac{2\pi h}{2L+1}, $$
and $\tilde{D}^{y,y'}_x$ is expressed according to (\ref{eq:D}) in terms of $\tilde{H}^y_x$. 

We should make several remarks. First, for sufficiently large $L$ one
can use the integral expression (\ref{eq:H}) for coefficients $H^y_x$
as an approximation. Unfortunately, such approximation cannot give an
accurate result for $|x|$ near $L$. When we are not interested in 
points near barriers, it does not lead to the problems. But in some
cases it can do. Therefore one should be careful working with finite membranes. 

Second, for absorbing and reflecting barriers there appear new unknown functions 
$P_{\pm L,y}(n)$ which can be called {\it barrier functions}. It means 
that one should write additional conditions to obtain a closed 
system of linear equations for near-boundary, ground and barrier functions. 
In order to obtain the condition corresponding to the barrier functions
$\phi _{\pm L,y}$, one can maintain simply $x=\pm L$ in (\ref{eq:supplem5}) 
taking different $y$. The problem is that there is an infinite number of barrier 
functions. A convenient solution is to 
consider two vertical barriers together with a horizontal barrier
(that is the usual case for a real physical membrane). In this case,
the number of barrier functions is just $\NN $. However, one should
not forget to use coefficients $H^y_x$ and $D^{y,y'}_x$ corrected by
introducing the horizontal barrier, see the previous subsection. 

Third, one can try to use a rough approximation for a very large $L$.
Indeed, if $|x|\ll L$, one has $D^{y,y'}_{x\pm L}\sim 1/L^2$, and
the barrier functions can be neglected. It simply means that one
uses the initial formula (\ref{eq:P}) to approximate the solution
with two remote barriers. 

Fourth, in the case of one horizontal and two vertical barriers 
there is a {\it finite} number of external points, i.e., one can
solve the system of discrete Laplacian equations (\ref{eq:Laplace_discrete})
properly, without using probabilistic techniques. In other words, the 
essential advantage of the present formalism is the possibility
to analyse the membranes with an infinite number of sites.

\vskip 5mm
\subsection{ Time-dependent distribution of hitting probabilities }                 \label{sec:time}

Previous results were related to the Laplacian equation, i.e.,
we considered distributions of hitting probabilities independent
of time. Generalizing the problem, we can find the distribution
of probabilities $P_{x,y}^{(t)}(n)$ to hit the $n$-th site of the 
membrane at $t$-th step (here $t$ is discrete time; for the 
Brownian motion it is continuous time) starting from an external 
point $(x,y)$. From the technical point of view, it is convenient
to consider the Laplacian transform of this distribution,
$$P_{x,y}^{(\lambda )}(n)=\sum\limits _{t=0}^{\infty }\lambda ^t
P_{x,y}^{(t)}(n).$$
Note that $\lambda =1$ corresponds to the previous time-independent
distribution. More generally, $P_{x,y}^{(\lambda )}(n)$
gives the distribution of hitting probabilities if 
particles have killing rate $(1-\lambda )$. 
For the planar case, one can apply the same technique as usual
to obtain the distribution
$$P_{x,y}^{(\lambda )}(n)=\int\limits _{-\pi }^{\pi }\frac{d\theta }{2\pi }
e^{i(x-n)\theta }\varphi ^y(\theta ;\lambda )=H^{y}_{x-n}(\lambda )$$
with function  
$$\varphi (\theta ;\lambda )=2/\lambda -\cos \theta -
\sqrt{(2/\lambda -\cos \theta )^2-1}.$$
This expression can be easily rewritten for the multi-dimensional case. 
Now we repeat the previous calculations for general membranes
to obtain exactly the same formula as (\ref{eq:P}),
$$P_{x,y}^{(\lambda )}(n)=\delta _{x,x_n}\delta _{y,y_n}+
D^{y,y_n+\delta y_n}_{x-x_n-\delta x_n}(\lambda )
-\sum\limits _{m=1}^M D^{y,y_m}_{x-x_m}(\lambda ) P_{[m]}^{(\lambda )}(n)
+\sum\limits _{g=1}^G H^y_{x-\bar{x}_g}(\lambda ) P_{\bar{x}_g,0}^{(\lambda )}(n),$$
where coefficients $D^{y,y'}_x(\lambda )$ are computed according to
formula (\ref{eq:D}) with $H^y_x(\lambda )$ instead of $H^y_x$. 
Imposing usual conditions for near-boundary and ground functions,
one can close these relations in order to find $P_{[m]}^{(\lambda )}(n)$
and $P_{\bar{x}_g,0}^{(\lambda )}(n)$. 
It is important that $\lambda $ appears in all expressions as a parameter. 

Now distribution $P_{x,y}^{(t)}(n)$ can be found as
$$P_{x,y}^{(t)}(n)=\frac{1}{t!}\biggl[\frac{d^t}{d\lambda ^t}
\bigl(P_{x,y}^{(\lambda )}(n)\bigr)\biggr]_{\lambda =0} .$$
Replacing function $\varphi (\theta ;\lambda )$ by 
$$\varphi (\theta _1,...,\theta _{d-1}; \lambda )=\left(\frac{d}{\lambda }-\sum\limits _{i=1}^{d-1}
\cos \theta _i\right)-\sqrt{\left(\frac{d}{\lambda }-\sum\limits _{i=1}^{d-1}
\cos \theta _i\right)^2-1} ,$$
we easily generalize this time-dependent distribution of hitting
probabilities to the multi-dimensional case.

\vskip 10mm
\section{ Numerical verifications }                                            \label{sec:numerical}

In this section we briefly present some numerical results 
to check the validity of the method.

\begin{figure}                   
\begin{center}
\includegraphics[width=125mm,height=5cm]{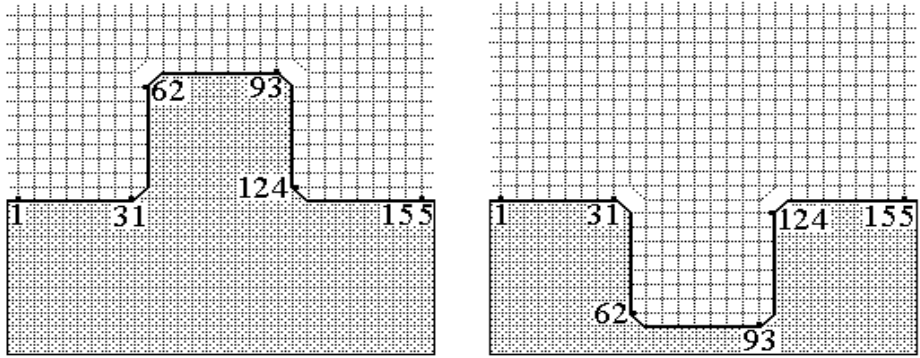}
\end{center}
\caption{ Two simple membranes : (a) convex membrane 
($N=31$, $\N =0$, $M=155$); \hskip 3mm (b) concave membrane 
($N=0$, $\N =31$, $M=155$). }                                                       \label{fig:2D}
\end{figure}

\begin{figure}                  
\begin{center}
\includegraphics[width=12cm,height=5cm]{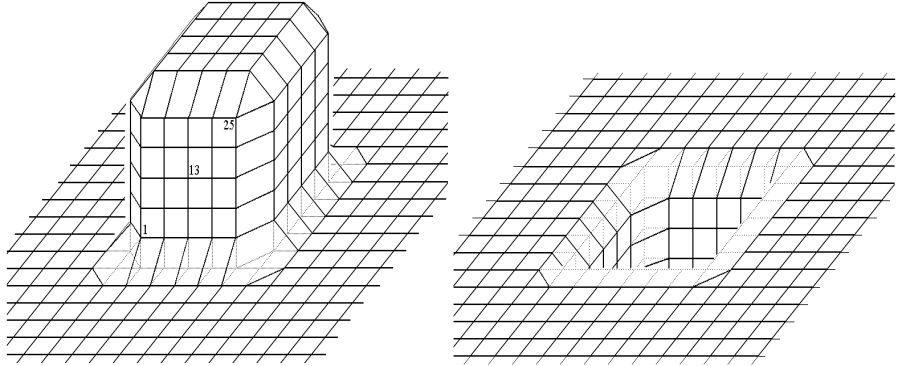}
\end{center}
\caption{ Three-dimensional membranes: (a) convex membrane ($N=6$, $\N =0$, $M=365$); 
(b) concave membrane ($N=0$, $\N =6$, $M=365$). There are $13$ facets $5\times 5$ 
(including $8$ facets on the ground level), and $8$ segments $5\times 1$.
{\small \it
Boundary points are enumerated by index $n$ following the rules: 
$1^{st}$ facet ($y=0$): $n=(|z|-1)\cdot 5+x$, \hskip 2mm
$2^{nd}$ facet ($x=6$): $n=(|z|-1)\cdot 5+y+25$,  \hskip 2mm  
$3^{td}$ facet ($y=6$): $n=(|z|-1)\cdot 5+x+50$,  \hskip 2mm 
$4^{th}$ facet ($x=0$): $n=(|z|-1)\cdot 5+y+75$,  \hskip 2mm 
$5^{th}$ facet ($|z|=6$): $n=(y-1)\cdot 5+x+100$, \hskip 2mm
$6^{th}-13^{th}$ facets ($z=0$) -- by the same way.}}                                \label{fig:3D}
\end{figure}

For a particular non-trivial membrane, the validity of expressions 
(\ref{eq:P}) and (\ref{eq:P_d}) can be verified by comparing them 
with the results obtained by numerical simulations (simple 
random walks on the lattice). We have taken two- and three-dimensional 
membranes represented on Fig.~\ref{fig:2D} and Fig.~\ref{fig:3D}.

\begin{figure}              
\begin{center}
\includegraphics[width=12cm,height=6cm]{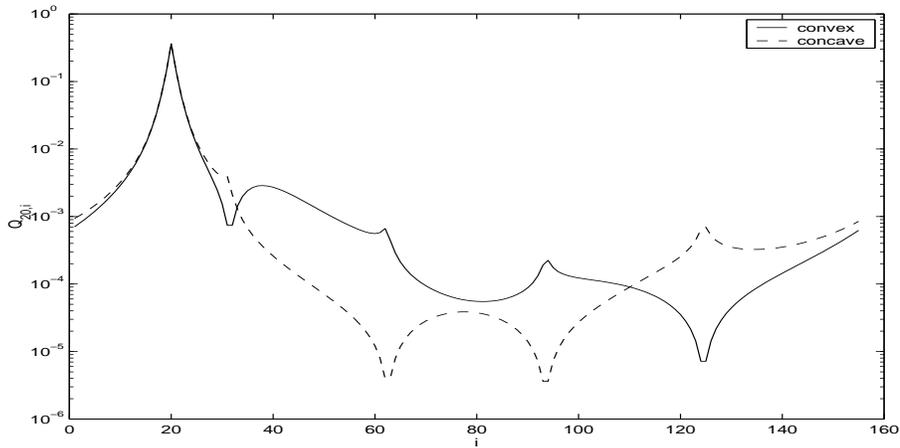}
\end{center}
\caption{ Distributions of hitting probabilities (in log-scale)
$Q_{k,n}$ with $k=20$: the probability of the first contact with 
points $(x_n,y_n)$ of membrane's boundary if started from the point 
$(x_{20},y_{20})$. }                                                            \label{fig:result2D}
\end{figure}

In the first case (Fig.~\ref{fig:2D}a), there are no ground functions, $G=0$, 
and there are $155$ near-boundary functions which can be calculated
using (\ref{eq:Q}). For the concave membrane (Fig.~\ref{fig:2D}b), there is $31$ 
ground functions and $155$ near-boundary functions.
On Fig.~\ref{fig:result2D} we present the distribution of hitting 
probabilities $P_{[20]}(n)$ (i.e., $Q_{20,n}$) obtained with the help 
of formula (\ref{eq:P}). The solid line corresponds to the convex membrane,
the dashed line -- to the concave membrane. These results are obtained
for the membranes placed in the finite box. In other words, we take the horizontal 
absorbing barrier on the distance $D=100$ from the ground level, and two vertical
barriers under the cyclic boundary conditions (see Section~\ref{sec:generalizations}).
Moreover, we account the corner points' corrections (see Appendices~\ref{sec:corner}). 

\begin{figure}             
\begin{center}
\includegraphics[width=12cm,height=6cm]{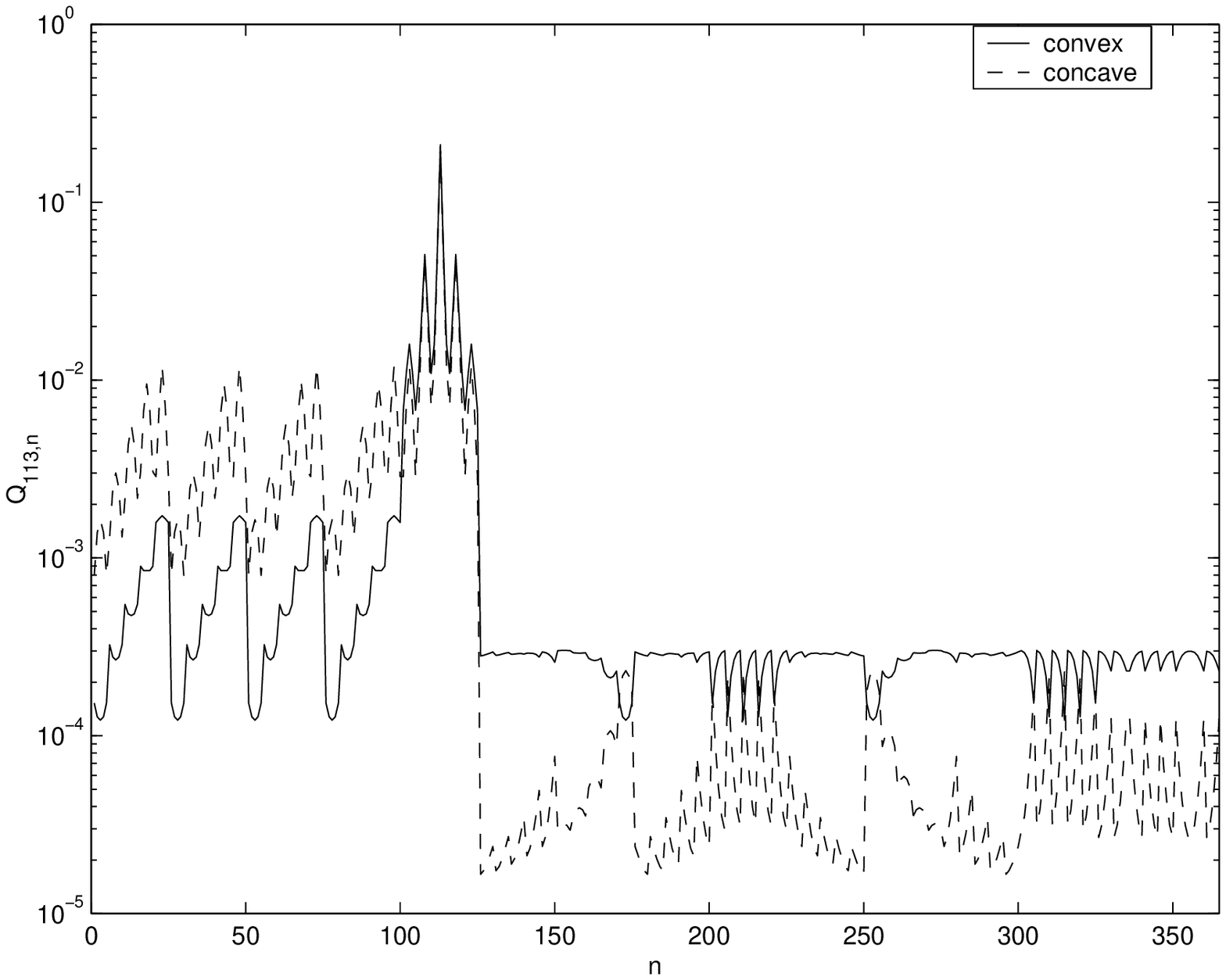}
\end{center}
\caption{ Distribution of hitting probabilities (in log-scale)
$Q_{k,n}$ with $k=113$ (the center of the fifth facet) 
for the three-dimensional convex and concave membranes.}                        \label{fig:result3D}
\end{figure}

Fig.~\ref{fig:result3D} shows the distribution of hitting probabilities $P_{[113]}(n)$
(i.e., $Q_{113,n}$) for three-dimensional membranes (see Fig.~\ref{fig:3D}): the solid 
line corresponds to the convex membrane, and the dashed line -- to the concave membrane.
These curves are obtained with the help of formula (\ref{eq:Q}), 
where we used (\ref{eq:H_d}) for $H^y_{\bx }$ instead of (\ref{eq:H}).
We calculate these distributions for the membranes placed in the finite box 
which is formed by the horizontal absorbing barrier at $D=15$,
and by four vertical barriers with cyclic boundary conditions
(see Section~\ref{sec:generalizations}). Also we took into account the corner 
points' corrections.

\vskip 2mm

In order to confirm that the method works, we performed the Monte-Carlo 
simulations. Starting from the near-boundary point $(x_k+\delta x_k,y_k+\delta y_k)$,
the particle walks at random on the lattice until it is absorbed either
on the membrane's boundary or on the horizontal barrier. 
Repeating this procedure, we obtain the numerical distribution of frequencies
of visits for all boundary points. Dividing by the total number of 
walks $N_t$, we calculate the approximate hitting probability $\tilde{Q}_{k,n}$ 
for each pair of sites $k$ and $n$. When $N_t$ goes to infinity,
the approximate values $\tilde{Q}_{k,n}$ tend to the exact values $Q_{k,n}$.
These simulations represent the well-known Bernoulli trials.
In order to estimate the accuracy of values $\tilde{Q}_{k,n}$ for 
the finite $N_t$, we use the central limit theorem. It gives the natural
measure $\Delta Q_{k,n}$ of deviations $Q_{k,n}-\tilde{Q}_{k,n}$, 
$$\Delta Q_{k,n}\sim \sqrt{Q_{k,n}}/\sqrt{N_t} ,$$
where $\sigma \approx \sqrt{Q_{k,n}}$ is the dispersion of the
Bernoulli trials.

Taking different values of $N_t$ and performing the Monte-Carlo simulations
for all possible sites $k$ and $n$, we compute the maximal deviation
of the approximate value $\tilde{Q}_{k,n}$ from the exact value $Q_{k,n}$,
$$\alpha =\max\limits _{k,n} \left\{ \frac{|Q_{k,n}-\tilde{Q}_{k,n}|}{\sqrt{Q_{k,n}}}\sqrt{N_t}\right \} .$$
Values of $\alpha $ for different membranes are given in Table 1.

\begin{center}
\begin{tabular}{| c | c | c | c |}  \hline
Membranes  & $10^5$  & $10^6$  & $10^7$  \\  \hline
2D convex  & $4.99$  & $4.02$  & $4.43$  \\
2D concave & $4.96$  & $4.91$  & $4.46$  \\
3D convex  & $5.33$  & $4.44$  & $4.46$  \\
3D concave & $6.33$  & $4.59$  & $4.46$  \\  \hline
\end{tabular}
\end{center}
{\em Table 1. Monte-Carlo verifications of the method's work. 
For $N_t\in \{ 10^5, 10^6, 10^7 \}$, the maximal deviation $\alpha $
is calculated.}

\vskip 1mm

These results confirm that formulae (\ref{eq:P}) and (\ref{eq:P_d}) work correctly.

\vskip 10mm
\section{ Discussion and conclusions }                                    \label{sec:conclusion}

In this paper, we defined a wide class of finite 
membranes which play an important role for the discrete Laplacian 
transfer problem. We based on the technique of characteristic functions 
developped in our previous paper. In particular, expressions (\ref{eq:conditions}) 
are used to close the recurrence relations (\ref{eq:recurrence}). 
We remind that expression $c_{N_u+1}\approx \varphi c_{N_u}$ was the 
central approximation of \cite{Grebenkov02}. In order to obtain this
approximate relation, we supposed that from a remote line $y=N_u$ 
(for large $N_u$) the membrane can be viewed as almost translationally 
invariant object (along the horizontal axis). Then the explicit solution 
of recurrence relations (\ref{eq:recurrence}) had been found. Finally, 
we took the limit $N_u\to \infty $ (and $N_l\to \infty $
for the lower half plane) to obtain expression (\ref{eq:c_y})
for coefficients $c_y(\theta ,\theta _h)$. It means that
influence of approximate relation (\ref{eq:conditions}) 
vanishes. In other words, the approximate relation
(\ref{eq:conditions}) is taken at infinity ($N_u\to \infty $), 
i.e., this condition has no influence on the solution.
Moreover, one can verify directly, without any approximation, 
that (\ref{eq:c_y}) is the exact solution of recurrence relations 
(\ref{eq:recurrence}). 

We stress that the essential result of our calculations 
is the form of solution (\ref{eq:P}) which expresses 
the distribution of hitting probabilities for a general
membrane, $P_{x,y}(n)$, in terms of corresponding
planar distribution, $H^y_x$. In other words, we 
explicitate how the membrane's geometry should change
the planar distribution. It means that having solved 
the planar problem, one can easily generalize its
solution for the case with a rather complex geometry.
For example, when the problem with horizontal
barrier for the planar case was solved, we just 
replaced old coefficients $H^y_x$ by the new ones.
The same ideas were used for the time-dependent 
distribution and for the multi-dimensional generalization.
Note that these motivations are frequently used
in the theory of conformal transforms. Indeed,
taking an appropriate conformal transform,
one can map the initial complex set into a simple
set (like a half plane or a disk), solve the 
simplified problem, and then reconstruct the solution
by the inverse conformal transform. Note also
that the conformal theory works only in two-dimensional
case while the present approach is valid for any
dimension of the lattice.

\vskip 2mm
\subsection{ Ground functions }                                              \label{sec:ground}

Let us briefly discuss the role of ground functions.
From the beginning, the upper and the lower half planes were considered
separately (for example, see two terms $\gamma ^{(y)}_{y'}$ and
$\gamma ^{(-y)}_{-y'}$ in equation (\ref{eq:c_y})).
The technique proposed in \cite{Grebenkov02} consists to {\it step down}
from the \hskip 1mm $N_u$-th and $(-N_l)$-th levels to the level zero.
It allows to express coefficients $c_y$ in terms of $\Delta c_y$.
But we cannot step down directly from $N_u$-th to $(-N_l)$-th level,
because there appears an infinity of near-boundary functions
due to $(\Delta \Phi _0)$, and we do not achieve our goal\footnote{ For the same
reason the compactness condition was imposed at the beginning.}.
In other words, we cannot ``pass'' through the level zero
in our treatment. On the other hand, a random walk started at the
upper half plane may hit the boundary point in the lower half
plane (if such point exists). This walk passes the level zero
through an external point, i.e., through a ground point. Therefore
we can say that ground functions ``connect'' the solutions for
the upper and lower half planes. 

Consider as example the concave membrane (Fig.~\ref{fig:2D}b). The probability
to hit the $n$-th site in the lower half plane if started from a point
$(x,y)$ in the upper half plane becomes
$$P_{x,y}(n)=\sum\limits _{g=1}^G H^y_{x-\bar{x}_g}P_{\bar{x}_g,0}(n) .$$
We simply took (\ref{eq:P}) under condition that $y>0$, $y_m<0$ for $m\in [1,M]$.
This formula has a clear probabilistic sense: to hit the $n$-th site,
the walk should reach one of the ground points (enumerated by $g$)
with probability\footnote{ We write $H^y_{x-\bar{x}_g}$ because in the upper 
half plane there is no ``perturbation'', the concave membrane completely lies 
in the lower half plane.} $H^y_{x-\bar{x}_g}$, and after that hit the
$n$-th site with probability $P_{\bar{x}_g,0}(n)$.

\vskip 2mm
\subsection{ Probabilistic sense of coefficients $D^{y,y'}_x$ }                   \label{sec:proba}

Now we clarify the nature of coefficients $D^{y,y'}_x$ which seem to be
artificial at the first sight. First of all, they {\it have no} purely 
probabilistic sense (for example, they can have values greater than $1$). 
To understand the sense of these coefficients, consider the simple problem:
to find the probability to hit a point $(x_0,y_0)$ if started from $(x,y)$
without touching the horizontal axis. The solution of this 
problem is given in Appendices 9.2,
\begin{equation}                                                                    \label{eq:P_one}
P_{[(x,y)\to (x_0,y_0)]}=\frac{D^{y,y_0}_{x-x_0}}{D^{y_0,y_0}_0 }.
\end{equation}
Using this solution, we can rewrite (\ref{eq:P}) for $(x,y)\notin \S $ as
\begin{equation}                                                                  \label{eq:Psimple}
P_{x,y}(n)=P_{[(x,y)\to (x_n+\delta x_n,y_n+\delta y_n)]}
D^{y_n+\delta y_n,y_n+\delta y_n}_0
\end{equation}
$$-\sum\limits _{m=1}^M P_{[(x,y)\to (x_m,y_m)]}D^{y_m,y_m}_0P_{[m]}(n)+
\sum\limits _{g=1}^G H^y_{x-\bar{x}_g}P_{\bar{x}_g,0}(n) .$$
Now this expression has more clear sense: we take all possible ways
to reach near-boundary point $(x_n+\delta x_n,y_n+\delta y_n)$ without
touching the horizontal axis, and then subtract all ways which go through
membrane's points. Coefficients $D^{y_m,y_m}_0$ and $D^{y_n+\delta y_n,y_n+\delta y_n}_0$
have combinatoric origin: they allow to account different possibilities
to go through the membrane. According to definition (\ref{eq:D}),
$$D^{y,y}_0=\sum\limits _{j=1}^y H^{2j-1}_0 ,$$
and for two-dimensional case the approximation
$$D^{y,y}_0\approx \frac{\ln y}{2\pi }+const \hskip 5mm (const\approx 0.3675)$$
gives $D^{y,y}_0$ with a sufficiently high accuracy.

Writing (\ref{eq:Psimple}) for the planar membrane, one has
$$P_{x,y}(n)=P_{[(x,y)\to (x_n,1)]}D^{1,1}_0 ,$$
i.e., a particle reaches near-boundary point $(x_n,1)$ with probability
$P_{[(x,y)\to (x_n,1)]}$, and then it hits corresponding boundary point $(x_n,0)$
with probability $D^{1,1}_0=H^1_0$ just as required.

\vskip 2mm
\subsection{ Coefficients $D^{y,y'}_{x-x'}$ like Green's functions }            \label{sec:Green}

Coefficients $D^{y,y'}_{x-x'}$ can be viewed from another standpoint.
A direct verification shows that
$$\L D^{y,y'}_{x-x'}=\delta _{x,x'}\delta _{y,y'} ,  \hskip 10mm  D^{y,0}_{x-x'}=0 ,$$
where the discrete Laplacian $\L $ operates on coordinates 
$(x,y)$ or $(x',y')$ (we remind that $D^{y,y'}_{x-x'}=D^{y',y}_{x'-x}$).
For example, in two-dimensional case we have
$$\L u_{x,y}=4u_{x,y}-u_{x+1,y}-u_{x-1,y}-u_{x,y+1}-u_{x,y-1}.$$
%
It means that coefficients $D^{y,y'}_{x-x'}$ can be treated as 
{\it Green's functions} of the discrete Dirichlet problem 
in the upper half plane. Consequently, the solution of
the general Dirichlet problem in the same (planar) geometry,
$$\L u_{x,y}=f_{x,y},  \hskip 10mm  u_{x,0}=0 ,$$
is given as
$$u_{x,y}=\sum\limits _{x',y'}D^{y,y'}_{x-x'} f_{x',y'} .$$
Unfortunately, this result has no direct connection
with our actual problem of general membranes. Nevertheless, 
we can use the following trick. What violates the using of
Green's function $D^{y,y'}_{x-x'}$? The answer is that there 
exist the boundary and internal points of the membrane with $y>0$.
Obviously, the Laplacian equation becomes invalid on these
points. However, we can correct this situation by introducing
a certain function $f_{x,y}(n)$ (here $n$ is a parameter). 
In other words, we impose the Laplacian equation artificially 
for the boundary and internal points.
The correction function is equal to $0$ for any external point. 
Moreover, $f_{x,y}$ equals to $0$ for almost all internal point
except the internal layer whose points are nearest to the membrane's 
boundary. For the boundary points we have
$$f_{x_m,y_m}(n)=4P_{x_m,y_m}(n)-P_{x_m-1,y_m}(n)-P_{x_m+1,y_m}(n)
-P_{x_m,y_m+1}(n)-P_{x_m,y_m-1}(n) ,$$
whereas for the nearest internal layer 
$$f_{x_m-\delta x_m,y_m-\delta y_m}(n)=-P_{x_m,y_m}(n) .$$
Now we can write the general solution in terms of Green's functions,
$$P_{x,y}(n)=\sum\limits _{m=1}^M \bigl[D^{y,y_m}_{x-x_m}f_{x_m,y_m}(n)
+D^{y,y_m-\delta y_m}_{x-x_m+\delta x_m}f_{x_m-\delta x_m,y_m-\delta y_m}(n) \bigr] .$$
Imposing the boundary condition $P_{x_m,y_m}(n)=\delta _{m,n}$ and
considering step by step four possible directions of the outer normal, 
one can demonstrate that this formula can be written as
$$P_{x,y}(n)=\delta _{x,x_n}\delta _{y,y_n}+D^{y,y_n+\delta y_n}_{x-x_n-\delta x_n}
-\sum\limits _{m=1}^M D^{y,y_m}_{x-x_m} P_{[m]}(n) .$$
Note that here we do not discuss neither contribution
of plane ``tails'', nor ground functions. One can easily 
complete these motivations in order to obtain formula (\ref{eq:P}).
This approach is also valid for multi-dimensional case.

\vskip 2mm
\subsection{ Comments on the Brownian self-transport operator }            \label{sec:comments}

The important application of the distribution (\ref{eq:P})
of hitting probabilities is the explicit construction of the Brownian
self-transport operator $Q$ which governes the Laplacian transfer.
We obtained the explicit formula (\ref{eq:Q}) for $Q$, and it opens the 
possibility of analytical researches in this field. Note that 
exact analytical expression of $Q$ even for a planar membrane
allowed to obtain the important characteristics of the 
Laplacian transfer, \cite{Grebenkov02_1}.

On the other hand, expression (\ref{eq:Q}) simplifies also
the numerical treatment of the problem. Usually one uses computer
simulations of random walks to calculate the matrix elements of $Q$.
In order to obtain the whole matrix $Q$ with a high accuracy, 
one should make enormous number of random walks. Moreover, 
one walk can be very long, especially on lattices with $d>2$.
On the contrary, working with our approach, one just
needs to manipulate in the framework of linear algebra. 
Once calculated, coefficients $H^y_{\bx }$ can be easily
used for any membrane. Therefore, one just needs to ``compose''
the matrix $D$ for a given membrane, and to inverse $(I+D)$.
The time required to make these operations depends only on 
the number of sites for which one calculates $Q$.

In order to adjust the formalism to a more realistic membrane,
one can introduce different barriers.  
As we described above, the horizontal barrier can be introduced 
by simple modification of coefficients $H^y_x$ while the problem 
with vectical barriers is more complicated.  

\vskip 1mm

We conclude that the present treatment opens a wide
field for further investigations.

\vskip 10mm
\section*{ Acknowledgement  }

The author would like to thank Professor B.Sapoval and Professor M.Filoche 
whose invaluable advice helped to develop the present approach and to construct  
the Brownian self-transport operator.

\vskip 10mm
\section{ Appendices }                                                      \label{sec:appendices}
\subsection{ Corner points' corrections }                                      \label{sec:corner}

In Section~\ref{sec:boundary}, we calculated the contributions of the boundary 
points supposing that these points are in general position. In other words, we did 
not consider the neighbourhood of the corners. Moreover, the procedure of corner points' 
removing distorts the lattice near corners (see Fig.~\ref{fig:corner}d), and we did 
not account  the influence of such distortion. For sufficiently regular membranes, 
these corrections should not change the distribution of hitting probabilities globally. 
However, the lattice distortions can perturb significantly the solution locally, 
in the vicinity of the corners. It can be important for certain problems. Here we 
briefly explain how solution (\ref{eq:P}) can be improved in order to account such 
local distortions. Note that these corrections carry a particular character, they depend 
on the way how the lattice has been distorted. For example, if one would not like to 
remove the corners, there is no correction.

\begin{figure}             
\begin{center}
\includegraphics[width=6cm,height=4cm]{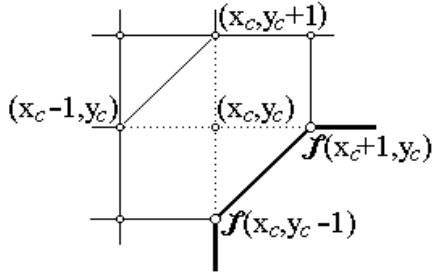}
\end{center}
\caption{ Corner points' corrections. }                                     \label{fig:correction}
\end{figure}

Why problems appear near the corners? Consider as example the corner point $(x_c,y_c)$
shown on Fig.~\ref{fig:correction}. 

We wrote the Laplacian equation (\ref{eq:relation}) for any
external point, in particular,  
$$4\phi _{x_c,y_c+1}-\phi _{x_c+1,y_c+1}-\phi _{x_c-1,y_c+1}-\phi _{x_c,y_c+2}-\phi _{x_c,y_c}=0 .$$
After that, we removed the corner point $(x_c,y_c)$ (thus one can take $\phi _{x_c,y_c}=0$),
and connected points $(x_c,y_c+1)$ and $(x_c-1,y_c)$ to each other by a liaison. It means 
that the true Laplacian equation is
$$4\phi _{x_c,y_c+1}-\phi _{x_c+1,y_c+1}-\phi _{x_c-1,y_c+1}-\phi _{x_c,y_c+2}-\phi _{x_c-1,y_c}=0 .$$
We prefer to work with the previous {\it general} form, and for these purposes 
we introduce here the correction term $\Delta \phi _{x_c,y_c+1}$ which is equal to 
$\phi _{x_c-1,y_c}$. The same concerns the point $(x_c-1,y_c)$. So, each external
corner gives two correction terms (in two-dimensional case). Note that these corrections 
have not been accounted above, they appear due to distortion of the lattice. 
We should also consider the Laplacian equation at point $(x_c,y_c)$ because
this {\it removed} point has been accounted in the general treatment,
$$4\phi _{x_c,y_c}-\phi _{x_c+1,y_c}-\phi _{x_c-1,y_c}-\phi _{x_c,y_c+1}-\phi _{x_c,y_c-1}=
\Delta \phi _{x_c,y_c} .$$
Above we supposed that the left hand side is equal to $0$ because the point $(x_c,y_c)$
was referred as external point. Now we calculate the corrections due to the lattice distortions, 
and for these purposes we write the correction term $\Delta \phi _{x_c,y_c}$. Taking all these terms, 
we obtain the whole contribution of the corner point $(x_c,y_c)$,
$$Z_c=-\bigl[e^{-ix_c\theta _h}\gamma ^{(y)}_{y_c}-e^{-ix_c\theta _h}\gamma ^{(y)}_{y_c+1}\bigr]
\phi _{x_c-1,y_c}-\bigl[e^{-ix_c\theta _h}\gamma ^{(y)}_{y_c}-e^{-i(x_c-1)\theta _h}
\gamma ^{(y)}_{y_c}\bigr]\phi _{x_c,y_c+1}$$
$$-e^{i\I (x_c,y_c-1)\theta }e^{-ix_c\theta _h}\gamma ^{(y)}_{y_c}-
e^{i\I (x_c+1,y_c)\theta }e^{-ix_c\theta _h}\gamma ^{(y)}_{y_c} .$$
Two last terms should be included into the contribution of boundary points with indices
$\I (x_c,y_c-1)$ and $\I (x_c+1,y_c)$ respectively (exactly these two terms are 
missed if one considers the boundary points near corners, see Section~\ref{sec:boundary}). 
Two first terms involve the unknown characteristic functions $\phi _{x_c-1,y_c}$ and 
$\phi _{x_c,y_c+1}$ that can be called {\it corner functions}. It means that in general case 
formula (\ref{eq:P}) becomes
$$P_{x,y}(n)=P_{x,y}^0(n)-\sum\limits _{corners}\biggl[\biggl(D^{y,y_c}_{x-x_c}-
D^{y,y_c+\delta y_c^{(2)}}_{x-x_c-\delta x_c^{(2)}}\biggr)P_{x_c+\delta x_c^{(1)},y_c+\delta y_c^{(1)}}(n)+$$
$$\biggl(D^{y,y_c}_{x-x_c}-D^{y,y_c+\delta y_c^{(1)}}_{x-x_c-\delta x_c^{(1)}}\biggr)
P_{x_c+\delta x_c^{(2)},y_c+\delta y_c^{(2)}}(n)\biggr] ,$$
where $P_{x,y}^0(n)$ denotes the right hand side of (\ref{eq:P}), and 
$(\delta x_c^{(1,2)},\delta y_c^{(1,2)})$ are two outer normals on the corner point $(x_c,y_c)$. 

The corner functions $P_{x_c+\delta x_c^{(1,2)},y_c+\delta y_c^{(1,2)}}(n)$ can be calculated 
by the same way that was used for near-boundary functions, i.e., we close the system of linear 
equations and solve it. The last formula gives a right distribution for all sites of the membrane, 
including the points near corners.

\vskip 5mm
\subsection{ One point problem }                                                \label{sec:one_point}

Consider the simple problem of finding the distribution of 
hitting probabilities $P_{x,y}(n)$ on the horizontal axis 
if there exists an absorbing point $(0,y_0)$ with $y_0>0$.
As consequence, we shall find the probability to hit the point
$(0,y_0)$ if started from $(x,y)$ without touching the
horizontal axis.

Here we cannot apply the technique of the planar case given that

-- there is no translational invariance ;

-- the Laplacian equation (\ref{eq:Laplace_discrete}) is invalid 
on the point $(0,y_0)$;

On the contrary, we can use formula (\ref{eq:c_y}). 
There exists the only correction $\Delta \Phi ^{(y_0)}$ for the
unique point $(0,y_0)$, i.e.,
$$\Delta c_{y_0}(\theta )=-\phi _{1,y_0}(\theta )-\phi _{-1,y_0}(\theta )
-\phi _{0,y_0+1}(\theta )-\phi _{0,y_0-1}(\theta )$$
(we omitted $4\phi _{0,y_0}$ because it equals to $0$). 
Using formula (\ref{eq:phi_xy}) and taking the limit $L\to \infty $, we obtain
$$\phi _{x,y}(\theta )=\int\limits _{-\pi }^{\pi }\frac{d\theta '}{2\pi }
e^{ix\theta '}\biggl(\varphi ^y(\theta ')c_0(\theta ,\theta ')+ 
\gamma ^{(y)}_{y_0}(\theta ')\Delta c_{y_0}(\theta )\biggr),$$
where 
$$c_0(\theta ,\theta ')=\sum\limits _{x=-\infty }^{\infty }e^{ix(\theta -\theta ')}
=2\pi \delta (\theta -\theta ') .$$
Applying the inverse Fourier transform, we obtain
\begin{equation}                                                                     \label{eq:Pone}
P_{x,y}(n)=H^{y}_{x-n}-D^{y,y_0}_x\biggl(P_{1,y_0}(n)+P_{-1,y_0}(n)+
P_{0,y_0+1}(n)+P_{0,y_0-1}(n)\biggr) .
\end{equation}
Looking on this relation, we say that there are four {\it near-boundary}
functions on four points around $(0,y_0)$. Taking corresponding $x$ and $y$
and summarizing four equations for near-boundary functions, 
we write the equation for their sum,
$$\sum\limits _{m=1}^4P_{[m]}(n)=\bigl(H^{y_0}_{1-n}+H^{y_0}_{-1-n}+H^{y_0+1}_{-n}
+H^{y_0-1}_{-n}\bigr)$$
$$-\bigl(D^{y_0,y_0}_1+D^{y_0,y_0}_{-1}+D^{y_0+1,y_0}_0+
D^{y_0-1,y_0}_0\bigr)\sum\limits _{m=1}^4P_{[m]}(n) .$$
Using the properties of coefficients $H^y_x$ and $D^{y,y'}_x$, we
can simplify this relation and express the sum of near-boundary functions,
$$\sum\limits _{m=1}^4 P_{[m]}(n)=\frac{H^{y_0}_{-n}}{D^{y_0,y_0}_0} .$$
Substituting this sum into (\ref{eq:Pone}), we have
$$P_{x,y}(n)=H^y_{x-n}-H^{y_0}_{-n}\frac{D^{y,y_0}_x}{D^{y_0,y_0}_0} .$$
So, the problem is solved. The probability to hit the point $(0,y_0)$
without touching the horizontal axis is
$$P^{(one)}_{x,y}=1-\sum\limits _{n=-\infty }^{\infty }P_{x,y}(n)=
\frac{D^{y,y_0}_x}{D^{y_0,y_0}_0} .$$

We hope that the consideration of this simple problem clarifies the sense of 
coefficients $D^{y,y'}_x$ and simplifies the understanding of general results.

\end{document}